\documentclass[11pt]{article}
\usepackage[dvips]{graphics}
\input{amssym.def}
\input{amssym.tex}

\title{\Large\bf A CRACK INDUCED BY A THIN RIGID INCLUSION\\
PARTLY DEBONDED FROM THE MATRIX}
\author{{\it by}  Y.A. ANTIPOV\\
antipov@math.lsu.edu\\ 
 {\it Department of Mathematics, Louisiana State University}\\
{\it Baton Rouge LA 70803, USA}\\
{\it and} S.M. MKHITARYAN\\
 {\it Department of Mechanics of Elastic and Viscoelastic Bodies}\\
{\it  National Academy of Sciences,
Yerevan 0019, Armenia} 
}

\newcommand{\bfm}[1]{\mbox{\boldmath ${#1}$}}

\newcommand{\beqa}{\begin{eqnarray}}
\newcommand{\eeqa}[1]{\label{#1}\end{eqnarray}}
\newcommand{\bequ}{\begin{equation}}
\newcommand{\eequ}[1]{\label{#1}\end{equation}}

\newcommand{\Md}{\partial}

\newcommand{\Ga}{\alpha}
\newcommand{\Gb}{\beta}
\newcommand{\Gd}{\delta}
\newcommand{\Ge}{\epsilon}
\newcommand{\Gve}{\varepsilon}
\newcommand{\Gf}{\phi}
\newcommand{\Gvf}{\varphi}
\newcommand{\Gg}{\gamma}
\newcommand{\Gc}{\chi}

\newcommand{\Gl}{\lambda}
\newcommand{\Gn}{\eta}

\newcommand{\Gs}{\sigma}

\newcommand{\Go}{\omega}
\newcommand{\Gx}{\xi}
\newcommand{\Gy}{\psi}

\newcommand{\GD}{\Delta}
\newcommand{\GF}{\Phi}
\newcommand{\GG}{\Gamma}
\newcommand{\GL}{\Lambda}

\newcommand{\GT}{\Theta}
\newcommand{\GS}{\Sigma}
\newcommand{\GU}{\Upsilon}
\newcommand{\GO}{\Omega}

\newcommand{\GY}{\Psi}

\newcommand{\BGo}{\bfm\omega}

\newcommand{\BGF}{\bfm\Phi}

\newcommand{\BGY}{\bfm\Psi}

\def\Bg{{\bf g}}

\def\BT{{\bf T}}

\newcommand{\beq}{\begin{equation}}
\newcommand{\eeq}{\end{equation}}
\newcommand{\barr}{\begin{eqnarray}}
\newcommand{\earr}{\end{eqnarray}}
\newcommand{\beqn}{\begin{equation*}}
\newcommand{\eeqn}{\end{equation*}}
\newcommand{\barrn}{\begin{eqnarray*}}
\newcommand{\earrn}{\end{eqnarray*}}

\newcommand{\fr}{\frac}

\newcommand{\supp}{\mathop{\rm supp}\nolimits}

\newcommand{\sgn}{\mathop{\rm sgn}\nolimits}

\newcommand{\I}{\mathop{\rm Im}\nolimits}

\newcommand{\const}{\mbox{const}}

\begin{document}
\maketitle


\begin{abstract}

The interaction of a thin rigid inclusion with a finite crack is studied. Two plane problems of elasticity 
are considered. The first one concerns the case when the upper side of the inclusion is completely debonded from the matrix, and the crack penetrates into the medium. In the second model, the upper side of the inclusion is partly separated from the matrix, that is the crack length $2a$ is less than $2b$,  the inclusion length. It is shown that both 
problems are governed by a singular integral equation of the same structure. Derivation of the closed-form
 solution of this integral equation  is the main result of the paper.  The solution is found by solving  the associated vector Riemann-Hilbert problem with the Chebotarev-Khrapkov matrix coefficient.
A feature of the method proposed is that the vector Riemann-Hilbert problem is set on a finite segment, while the 
original Khrapkov method of matrix factorization is developed for a closed contour. 
In the case, when the crack and inclusion lengths are the same, the solution is derived by passing to the limit $b/a\to 1$. It is demonstrated that the limiting case $a=b$ is unstable, and  when $a<b$, and the crack  tips approach the inclusion ends, the crack tends to accelerate in order to penetrate into the matrix.

\end{abstract}


\setcounter{equation}{0}

\section{Introduction} 

Different aspects of interaction of cracks with inclusions have been studied in detail by many investigators.
The fundamental two-dimensional elastic model problem for a finite set of slits lying on the real axis when the traction components
and the displacements are prescribed in the upper and lower surfaces, respectively, was solved in  {\bf(\ref{she})}
by the method of singular integral equations and in  {\bf(\ref{mus})} by the method of complex potentials. 
Some generalizations of this setting and the method of the vector Riemann-Hilbert problem
were proposed and the interaction of a semi-infinite inclusion and a finite cut
was analyzed in  {\bf(\ref{cher})}.   The model {\bf(\ref{she})} for an inclusion whose upper side is completely separated from the matrix was
generalized in {\bf(\ref{pop})} for the case when an inclusion is in the interface between two half-planes
with different elastic constants.
The method of complex potentials and the scalar Riemann-Hilbert problem on a 
hyperelliptic surface was advanced  {\bf(\ref{zve})} for the mixed boundary value problem governing a system of collinear 
cuts when the points of change of boundary conditions are not necessarily the endpoints of the cuts.  
This technique was applied  {\bf(\ref{sil})} for the solution of the plane problem for $n$ inclusions
whose upper sides are completely separated from the matrix, and  slip lines  emerge from the crack tips.
The case of a rigid inclusion $-a<x<a, y=0$  with a crack $-a<x<a$ on  the inclusion upper side 
open in the interval $-b<x<b$ and contacting with the inclusion when $b<|x|<a$: $\Gs_{12}=0$, $u_2=0$
was treated in {\bf(\ref{ant})}.
A closed-form solution to the plane problem on a rigid inclusion $-b<x<b$
located on the lower side of a crack $-a<x<a, y=0$ ($b<a$) when the tangential  traction component vanishes
on the segment $-b<x<b, y=0^-$) (frictionless contact of the inclusion and the matrix) was found in {\bf(\ref{bar})}. Some other relevant works
include {\bf(\ref{sel})} and  {\bf(\ref{hwu})}. In the former paper, the problem of the symmetric indentation of a penny-shaped crack by a smoothly
embedded rigid circular thin disc  was treated by the method of triple integral equations. The second work
examines
models  of a crack inside, outside, penetrating or lying along the interface of an anisotropic elliptical inclusion 
and presents numerical solutions of the governing singular integral equations.

The main goal of this paper was to develop a method of 
 the vector Riemann-Hilbert on a segment for singular
  integral equation governing the following two model problems. 
Problem  1 concerns a rigid inclusion $-b<x<b$, $y=0$, whose lower side is bonded to
the matrix, while the upper side is completely separated from the elastic medium, and the crack
formed on the upper surface penetrates into the medium and occupies the segment $-a<x<a$, $y=0$ 
($a>b$). The latter is on a rigid inclusion $-b<x<b$, whose upper part is  debonded on the segment
$-a<x<a$,  and $a<b$.  
In particular, on passing to the limit $b/a\to 1$ we aim to examine the transition of the solution
to the first problem with the square root singularities at the crack tips to the solution
for the case $a=b$ that oscillates and has a stronger singularity of order $-3/4$ at the endpoints.

In Section 2, we formulate the first problem as a discontinuous boundary value problem for the biharmonic 
operator. It is reduced to a singular integral equation in Section 3. Then, in Section 4, we transform
the integral equation to a vector Riemann-Hilbert problem on a segment with the Chebotarev-Khrapkov
matrix coefficient {\bf(\ref{che})}, {\bf(\ref{khr})} and construct the solution of the Riemann-Hilbert problem.
Based on this solution we derive the closed-form solution of the singular integral equation in Section 5.
In the next section, we analyze the singularities of the solution, discuss the numerical results obtained
and find the solution in the limiting case $b/a\to 0$.
In Section 7, we show that Problem 2 is governed by a singular integral equation
whose kernel is the same as the one for the first problem. What is different is the right-hand side,
the additional condition and the meaning of the parameters and the unknown function. 
Also, in this section we implement the passage to the limit $a\to b$, and show that
in a particular case the limiting solution coincides with the solution available in the literature
{\bf(\ref{mus})}.

\setcounter{equation}{0}

\section{Formulation}\label{form}

\begin{figure}[t]
\centerline{
\scalebox{0.6}{\includegraphics{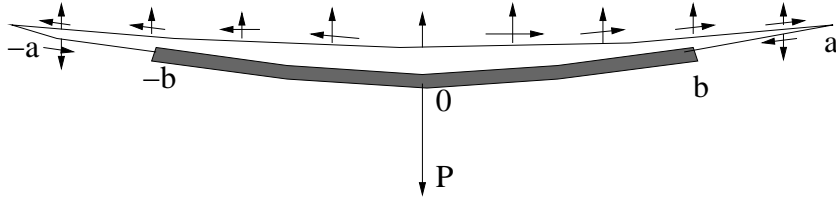}}
}
\caption{The geometry of Model 1.}
\label{fig1}
\end{figure}

The problem to be addressed is a two-dimensional one for a thin rigid inclusion 
whose profile is described by a function $y=h(x)$, $-b<x<b$. The function $h(x)$ is even,
$h(\pm b)=0$, and the curvature of $h(x)$, $\Ge$,  is small, $\Ge<<1$. The medium is assumed to be elastic, uniform, infinite and
to be loaded in a way that would generate, in the absence  of the inclusion, a stress field
$\Gs_{ij}^\circ(x,y)$ symmetric with respect to both axes, $x$ and $y$. The Poisson ratio and the Young modulus of the medium are $\nu$ and $E$, respectively, and the conditions of plane stress are considered. 
It is  assumed that along the upper side of the inclusion (Fig.1)
there is a crack that spreads not only over the whole upper surface of the inclusion, but also penetrates
into the matrix. We aim to analyze the stress concentration in the vicinity of the tips of the inclusion and the crack.

Owing to small deviations of the function $h(x)$ from the $x$-axis, it is conventional 
in linear elasticity to write the boundary conditions on the line $y=0$, not the actual curve $y=h(x)$.
They read
$$
\Gs_{12}+i\Gs_{22}=-\Gs_{12}^\circ-i\Gs_{22}^\circ, \quad |x|<a, \quad y=0^+,
$$$$
\Gs_{12}+i\Gs_{22}=-\Gs_{12}^\circ-i\Gs_{22}^\circ, \quad b<|x|<a, \quad y=0^-,
$$
\beq
\fr{\Md u_1}{\Md x}+i\fr{\Md u_2}{\Md x}=ih'(x)- w^\circ(x),  \quad |x|<b, \quad y=0^-.
\label{2.1}
\eeq
Here, $(u_1,u_2)$ is the displacement vector, $w^\circ(x)=\fr{\Md}{\Md x}(u_1^\circ+iu_2^\circ)(x,0)$, and $(u^\circ_1,u^\circ_2)$ is the displacement vector associated with the stress field $\Gs_{ij}^\circ$.

It will be convenient to introduce 
the functions
\beq
\Gf_j(x)=E\left[\fr{\Md u_j}{\Md x}\right](x), \quad \psi_j(x)=[\Gs_{j2}](x), \quad j=1,2,
\label{2.2}
\eeq
where
\beq
[f](x)=f(x,0^+)-f(x,0^-).
\label{2.3}
\eeq
Since the displacements are discontinuous across the whole crack surface and the traction
components are discontinuous across the inclusion surface and continuous across the segments
$-a<x<-b$ and $b<x<a$, we have $\supp \Gf_j\subset[-a,a]$ and $\supp\psi_j\subset[-b,b]$, $j=1,2$.

Let $U(x,y)$ be the Airy function of the problem. Then
$$
\Gs_{12}=-\fr{\Md^2 U}{\Md x\Md y}, \quad \Gs_{22}=\fr{\Md^2 U}{\Md x^2}, 
$$
$$
\fr{\Md u_1}{\Md x}=\fr{1}{E}\left(\fr{\Md^2 U}{\Md y^2}-\nu\fr{\Md^2 U}{\Md x^2}\right),
\quad 
\fr{\Md u_2}{\Md y}=\fr{1}{E}\left(\fr{\Md^2 U}{\Md x^2}-\nu\fr{\Md^2 U}{\Md y^2}\right),
$$
\beq
\fr{\Md u_1}{\Md y}+\fr{\Md u_2}{\Md x}=-\fr{2(1+\nu)}{E}\fr{\Md^2 U}{\Md x\Md y},
\label{2.4}
\eeq
and the model problem is equivalent to the following discontinuous boundary value problem for the
biharmonic operator:
$$
\GD^2 U(x,y)=0, \quad (x,y)\in {\Bbb R}^2\setminus\{-a<x<a, y=0\},
$$$$
[U](x)=\tilde\psi_2(x), \quad \left[\fr{\Md U}{\Md y}\right](x)=-\hat\psi_1(x),
$$
\beq
\left[\fr{\Md^2 U}{\Md y^2}\right](x)= \Gf_1(x)+\nu\psi_2(x),\quad 
\left[\fr{\Md^3 U}{\Md y^3}\right](x)= -\Gf_2'(x)+(\nu+2)\psi_1'(x),
\label{2.5}
\eeq
subject to the boundary conditions (\ref{2.1}), the inclusion equilibrium conditions
\beq
\int_{-b}^b\psi_1(x)dx=0, \quad \int_{-b}^b\psi_2(x)dx=\GS+P,
\label{2.6}
\eeq
and the crack closedness  conditions
\beq
\int_{-a}^a\Gf_1(x)dx=0, \quad \int_{-a}^a\Gf_2(x)dx=0.
\label{2.6'}
\eeq
Here, $\hat\psi_1'(x)=\psi_1(x)$ and $\tilde\psi_2''(x)=\psi_2(x)$, 
\beq
P=-\int_{-b}^b\Gs_{22}(x,0^-)dx, \quad \GS=-\int_{-b}^b\Gs_{22}^\circ(x,0^+)dx,\quad j=1,2,
\label{2.6''}
\eeq
$P$ is the magnitude of the total normal force applied to the inclusion central point and directed downwards.

\setcounter{equation}{0}
  
\section{Singular integral equation}\label{SIE}

Aiming to reduce the  discontinuous boundary value problem to an integral equation we apply first the 
Fourier transform across the discontinuity line 
\beq
\hat U_\Ga(x)=\int_{-\infty}^\infty U(x,y)e^{i\Ga y}dy
\label{3.1}
\eeq
and deduce the following ordinary differential equation for the Fourier transform $\hat U_\Ga(x)$:
$$
\left(\fr{d^4}{dx^4}-2\Ga^2\fr{d^2}{dx^2}+\Ga^4\right)\hat U_\Ga(x)
$$
\beq
=-i\Ga\Gf_1(x)-\Gf_2'(x)+\nu\psi_1'(x)+\Ga^2\hat\psi_1(x)-i\Ga(2+\nu)\psi_2(x)+i\Ga^3\tilde\psi_2(x), 
\quad -\infty<x<\infty.
\label{3.2}
\eeq
We next utilize the fundamental function of the differential operator in equation (\ref{3.2})
\beq
\fr{1+|\Ga||x-\Gx|}{4|\Ga|^3}e^{-|\Ga||x-\Gx|}
\label{3.3}
\eeq
and integration by parts to derive
$$
\hat U_\Ga(x)=\fr{1}{4\Ga^2}\int_{-a}^a\{-i\sgn\Ga (1+|\Ga||x-\xi|)\Gf_1(\xi)+|\Ga|(x-\xi)\Gf_2(\xi)
$$
\beq
-[2\sgn(x-\xi)+(1+\nu)
|\Ga|(x-\xi)]\psi_1(\xi)+i\sgn\Ga[1-\nu-(1+\nu)|\Ga||x-\xi|]\psi_2(\xi)\}e^{-|\Ga||x-\xi|}d\xi.
\label{3.4}
\eeq
On inverting the Fourier transform, using formulas (\ref{2.4}) and evaluating the corresponding integrals, 
we obtain the integral representations of the stresses 
$$
\Gs_{12}(x,y)=\int_{-a}^a
\{-M_1(x-\xi,y)\Gf_1(\xi)+[N_2(x-\xi,y)-M_2(x-\xi,y)]\Gf_2(\xi)\}d\xi
$$$$
+\int_{-b}^b
\{[(1-\nu)N_2(x-\xi,y)+(1+\nu)M_2(x-\xi,y)]\psi_1(\xi)+[N_1(x-\xi,y)-(1+\nu)M_1(x-\xi,y)]\psi_2(\xi)\}d\xi,
$$
$$
\Gs_{22}(x,y)=\int_{-a}^a
\{[N_2(x-\xi,y)-M_2(x-\xi,y)]\Gf_1(\xi)
-[N_1(x-\xi,y)-M_1(x-\xi,y)]\Gf_2(\xi)\}d\xi
$$
\beq
+\int_{-b}^b
\{[\nu N_1(x-\xi,y)-(1+\nu)M_1(x-\xi,y)]\psi_1(\xi)+[(3+\nu)N_2(x-\xi,y)-(1+\nu)M_2(x-\xi,y)]\psi_2(\xi)\}d\xi,
\label{3.5}
\eeq
and the displacement derivatives
$$
E\fr{\Md u_1}{\Md x}(x,y)=\int_{-a}^a
\{[(1-\nu)N_2(x-\xi,y)+(1+\nu)M_2(x-\xi,y)]\Gf_1(\xi)+[\nu N_1(x-\xi,y)
$$$$
-(1+\nu)M_1(x-\xi,y)]\Gf_2(\xi)\}d\xi
+\int_{-b}^b\{
[(1-\nu^2) N_1(x-\xi,y)+(1+\nu)^2M_1(x-\xi,y)]\psi_1(\xi)
$$$$
-(1+\nu)^2[N_2(x-\xi,y)-M_2(x-\xi,y)]\psi_2(\xi)
\}d\xi,
$$
$$
E\fr{\Md u_2}{\Md x}(x,y)=\int_{-a}^a
\{
[N_1(x-\xi,y)-(1+\nu)M_1(x-\xi,y)]\Gf_1(\xi)
+[(3+\nu)N_2(x-\xi,y)
$$$$
-(1+\nu)M_2(x-\xi,y)]\Gf_2(\xi)\}d\xi
+\int_{-b}^b\{
(1+\nu)^2[-N_2(x-\xi,y)+M_2(x-\xi,y)]\psi_1(\xi)
$$
\beq
+(1+\nu)[2 N_1(x-\xi,y)
-(1+\nu)M_1(x-\xi,y)]\psi_2(\xi)\}d\xi
\label{3.6}
\eeq
in terms of the 
jumps introduced in (\ref{2.2}). Here,
$$
N_1(t,y)=\fr{t}{2\pi(t^2+y^2)}, \quad N_2(t,y)=\fr{y}{4\pi(t^2+y^2)},
$$
\beq
M_1(t,y)=\fr{t(t^2-y^2)}{4\pi(t^2+y^2)^2}, \quad M_2(t,y)=\fr{t^2y}{2\pi(t^2+y^2)^2}.
\label{3.7}
\eeq
To satisfy the boundary conditions (\ref{2.1}), we need the boundary values 
$\Gs_{j2}(x,0^\pm)$ and $\fr{\Md}{\Md x}u_j(x,0^\pm)$, $j=1,2$. In view of the properties
of the Dirac delta-function
\beq
\lim_{y\to 0^\pm}\fr{y}{\pi(t^2+y^2)}=\pm\Gd(t), \quad 
\lim_{y\to 0^\pm}\fr{yt^2}{\pi(t^2+y^2)^2}=\pm\fr12\Gd(t),
\label{3.8}
\eeq
the integral representations (\ref{3.5}) and (\ref{3.6}) yield
$$
\Gs_{12}(x,0^\pm)=-\fr{1}{4\pi}\int_{-a}^a
\fr{\Gf_1(\Gx)d\Gx}{x-\Gx}\pm\fr12\psi_1(x)+
\fr{1-\nu}{4\pi}\int_{-b}^b\fr{\psi_2(\Gx)d\Gx}{x-\Gx},
$$
$$
\Gs_{22}(x,0^\pm)=-\fr{1}{4\pi}\int_{-a}^a
\fr{\Gf_2(\Gx)d\Gx}{x-\Gx}-
\fr{1-\nu}{4\pi}\int_{-b}^b\fr{\psi_1(\Gx)d\Gx}{x-\Gx}\pm\fr12\psi_2(x),
$$
$$
E\fr{\Md u_1}{\Md x}(x,0^\pm)=\pm\fr12\Gf_1(x)-\fr{1-\nu}{4\pi}\int_{-a}^a
\fr{\Gf_2(\Gx)d\Gx}{x-\Gx}+
\fr{\nu_1}{4\pi}\int_{-b}^b\fr{\psi_1(\Gx)d\Gx}{x-\Gx},
$$
\beq
E\fr{\Md u_2}{\Md x}(x,0^\pm)=\fr{1-\nu}{4\pi}\int_{-a}^a
\fr{\Gf_1(\Gx)d\Gx}{x-\Gx}\pm\fr12\Gf_2(x)+
\fr{\nu_1}{4\pi}\int_{-b}^b\fr{\psi_2(\Gx)d\Gx}{x-\Gx},
\label{3.9}
\eeq
where 
\beq
\nu_1=(3-\nu)(1+\nu).
\label{3.9'}
\eeq
It is convenient to introduce the following complex-valued functions:
$$
\Gf(x)=\Gf_1(x)+i\Gf_2(x),\quad \Gy(x)=\Gy_1(x)+i\Gy_2(x),
$$
\beq
\Gs_\pm(x)=\Gs_{12}(x,0^\pm)+i\Gs_{22}(x,0^\pm), \quad 
w_\pm(x)=\fr{\Md u_1}{\Md x}(x,0^\pm)+i\fr{\Md u_2}{\Md x}(x,0^\pm).
\label{3.10}
\eeq
Equations (\ref{3.9}), after rearrangement, become
$$
\Gs_\pm(x)=\pm\fr12\psi(x)+
\fr{i(1-\nu)}{4\pi}\int_{-b}^b\fr{\psi(\Gx)d\Gx}{\Gx-x}+\fr{1}{4\pi}\int_{-a}^a
\fr{\Gf(\Gx)d\Gx}{\Gx-x},
$$
\beq
Ew_\pm(x)=\pm\fr12\Gf(x)-\fr{i(1-\nu)}{4\pi}\int_{-a}^a
\fr{\Gf(\Gx)d\Gx}{\Gx-x}-
\fr{\nu_1}{4\pi}\int_{-b}^b\fr{\psi(\Gx)d\Gx}{\Gx-x}.
\label{3.11}
\eeq
We now proceed with the boundary conditions on the crack sides and the lower side of the inclusion.
From (\ref{3.11}), we have 
$$
2\Gs_+(x)=\fr{1}{2\pi}\int_{-a}^a
\fr{\Gf(\Gx)d\Gx}{\Gx-x}+\psi(x)+\fr{i(1-\nu)}{2\pi}\int_{-b}^b\fr{\psi(\Gx)d\Gx}{\Gx-x},
\quad -a<x<a,
$$
\beq
2Ew_-(x)=-\Gf(x)-i(1-\nu)[2\Gs_+(x)-\psi(x)]-\fr{2}{\pi}\int_{-b}^b\fr{\psi(\Gx)d\Gx}{\Gx-x},
\quad -b<x<b.
\label{3.12}
\eeq
On inverting the singular operator in the space of functions having integrable singularities
at the endpoints $\pm a$ {\bf(\ref{gak})} we express the function $\Gf(x)$ from the first equation in (\ref{3.12})
\beq
\Gf(x)=-\fr{2}{\pi\sqrt{a^2-x^2}}\left\{\int_{-a}^a \fr{\sqrt{a^2-\Gx^2}}{\Gx-x}
[2\Gs_+(\Gx)
-\psi(\Gx)]d\Gx+C\right\}+i(1-\nu)\GT(x),
\label{3.13}
\eeq
where $C$ is an arbitrary constant and 
\beq
\GT(x)=\fr{1}{\pi\sqrt{a^2-x^2}}\int_{-a}^a \fr{\sqrt{a^2-\Gx^2}}{\Gx-x}
\fr{1}{\pi}\int_{-a}^a\fr{\psi(\Gn)d\Gn}{\Gn-\Gx}.
\label{3.14}
\eeq
It turns out that $C$ has to be zero. Indeed, from the crack closedness conditions (\ref{2.6'})  and the  definition (\ref{3.10}) of the function $\Gf(x)$ we have
\beq
\int_{-a}^a\Gf(x)dx=0.
\label{3.15}
\eeq
On the other hand, 
\beq
\int_{-a}^a \fr{dx}{\sqrt{a^2-x^2}(x-\Gx)}=0, \quad -a<\Gx<a,
\label{3.16}
\eeq
where the integral in the left-hand side is understood in the sense of the Cauchy principal value. Therefore, from the representation (\ref{3.13}) of $\Gf(x)$  it follows that $C=0$. Our next step is to
simplify the relation (\ref{3.14}) for the function $\GT(x)$. On employing the Poincar\'e-Bertrand formula  {\bf(\ref{gak})}
\beq
\fr{1}{\pi}\int_{-a}^a \fr{d\Gx}{\Gx-x}\fr{1}{\pi}\int_{-a}^a \fr{\Gvf(\Gx,\Gn)d\Gn}{\Gn-\Gx}=
-\Gvf(x,x)+\fr{1}{\pi}\int_{-a}^a \fr{d\Gn}{\pi}\int_{-a}^a \fr{\Gvf(\Gx,\Gn)d\Gx}
{(\Gx-x)(\Gn-\Gx)},
\label{3.17}
\eeq
evaluating the integral
\beq
\fr{1}{\pi}\int_{-a}^a\fr{\sqrt{a^2-\Gx^2}d\Gx}{(\Gx-x)(\Gn-\Gx)}=1, 
\quad -a<x<a, \quad -a<\Gn<a,
\label{3.18}
\eeq
and taking into account the equilibrium conditions (\ref{2.6}) and the formula for the function
$\Gy(x)$ in (\ref{3.10}) we arrive ultimately at the following simple relation:
\beq
\GT(x)=-\psi(x)+\fr{iP_*}{\pi\sqrt{a^2-x^2}}, \quad -a<x<a,
\label{3.19}
\eeq
where $P_*=P+\GS$.
Now we replace the function $\Gf(x)$ in the second equation (\ref{3.12}) by
its representation 
$$
\Gf(x)=\fr{1}{\pi\sqrt{a^2-x^2}}\left[2\int_{-b}^b \fr{\sqrt{a^2-\Gx^2}}{\Gx-x}\psi(\Gx)d\Gx
+4\int_{-a}^a \fr{\sqrt{a^2-\Gx^2}}{\Gx-x}\Gs^\circ(\Gx)d\Gx-(1-\nu)P_*\right]
$$
\beq
-i(1-\nu)\Gy(x),\quad -a<x<a.
\label{3.19'}
\eeq
following from 
(\ref{3.13}) and obtain a governing singular integral equation
for the function $\Gy(x)$. It reads
\beq
\fr{1}{\pi}\int_{-k}^k\left(1+\sqrt{\fr{1-\tau^2}{1-t^2}}
\right)\fr{\Gy(a\tau)d\tau}{\tau-t}
-i(1-\nu)\Gy(at)=g(t), \quad -k<t<k.
\label{3.20}
\eeq
Here, $k=b/a<1$,
\beq
g(t)=-iEh'(at)+Ew^\circ(at)+i(1-\nu)\Gs^\circ(at)-\fr{2}{\pi\sqrt{1-t^2}}\int_{-1}^1
\fr{\sqrt{1-\tau^2}}{\tau-t}\Gs^\circ(a\tau)d\tau+\fr{(1-\nu)P_*}{2\pi a\sqrt{1-t^2}},
\label{3.21}
\eeq
$\Gs^\circ(x)=\Gs_{12}^\circ(x,0)+i\Gs_{22}^\circ(x,0)$, and due to (\ref{2.6})
the function $\psi(a\tau)$
has to satisfy the additional condition
\beq
\int_{-k}^k\psi(at)dt=\fr{iP_*}{a}.
\label{3.22}
\eeq

\setcounter{equation}{0}
  
\section{Vector Riemann-Hilbert problem and its solution}\label{RH} 

In this section we aim to rewrite the singular integral equation (\ref{3.20}) as a vector Riemann-Hilbert problem 
and derive its closed-form solution.
Introduce two Cauchy integrals
\beq
\GF_1(z)=\fr{1}{2\pi i}\int_{-k}^k\fr{\psi(a\tau)d\tau}{\tau-z}, \quad 
\GF_2(z)=\fr{1}{2\pi i}\int_{-k}^k\fr{\sqrt{1-\tau^2}\psi(a\tau)d\tau}{\tau-z}, 
\label{4.1}
\eeq
analytic in the whole complex plane except at the segment $[-k,k]$. By the Sokhotski-Plemelj formulas 
their limiting values $\GF_j^\pm(t)=\GF_j(t\pm i0)$, $-k<t<k$, satisfy the relations 
 \beq
\GF_1^+(t)+\GF_1^-(t)=\fr{1}{\pi i}\int_{-k}^k
\fr{\psi(a\tau)d\tau}{\tau-t},\quad
\GF_2^+(t)+\GF_2^-(t)=\fr{1}{\pi i}\int_{-k}^k
\fr{\sqrt{1-\tau^2}\psi(a\tau)d\tau}{\tau-t},
\label{4.2}
\eeq
and 
\beq
\GF_1^+(t)-\GF_1^-(t)=\psi(at), \quad \GF_2^+(t)-\GF_2^-(t)=\sqrt{1-t^2}\psi(at).
\label{4.3}
\eeq
On substituting these relations into equation (\ref{3.20}) we obtain that the limiting values of the functions
(\ref{4.1}) are connected by
$$
\nu\GF_1^+(t)+(2-\nu)\GF_1^-(t)+\fr{\GF_2^+(t)+\GF_2^-(t)}{\sqrt{1-t^2}}=-ig(t),
$$
\beq
\GF_2^+(t)-\GF_2^-(t)=\sqrt{1-t^2}[\GF_1^+(t)-\GF_1^-(t)], \quad -k<t<k.
\label{4.4}
\eeq
This, after rearrangement, can be written as the following vector Riemann-Hilbert boundary value 
problem on a segment:
\beq
\BGF^+(t)=G(t)\BGF^-(t)-\fr{i}{\nu+1}\BT (t)g(t), \quad -k<t<k,
\label{4.5}
\eeq
where
$$
G(t)=\fr{1}{\nu+1}\left(
\begin{array}{cc}
\nu-1 \; & \; -\fr{2}{\sqrt{1-t^2}}\\
-2\sqrt{1-t^2} \; & \; \nu-1 \\
\end{array}\right), 
 $$
 \beq
 \BGF(z)=\left(
\begin{array}{c}
\GF_1(z)\\
\GF_2(z) \\
\end{array}\right), \quad  \BT(t)=\left(
\begin{array}{c}
1\\
\sqrt{1-t^2} \\
\end{array}\right).
\label{4.6}
\eeq
The vector $\GF(z)$ has a simple zero at the infinite point, and due to the
condition (\ref{3.22}) and the first relation (\ref{4.1}) the function
$\GF_1(z)$ behaves at infinity as
\beq
\GF_1(z)\sim-\fr{P_*}{2\pi az}, \quad z\to\infty.
\label{4.6'}
\eeq 
The matrix coefficient to be factorized has the Chebotarev-Khrapkov structure  {\bf(\ref{che})},  {\bf(\ref{khr})},  and its eigenvalues
are constants, $\Gl_1=-(3-\nu)(1+\nu)^{-1}$, $\Gl_2=1$. We point out that the classical scheme  {\bf(\ref{khr})} 
is designed and has been employed in the literature so  far for the case of a closed contour. 
 In the case of a segment, although the structure of the Wiener-Hopf factors is preserved,
additional efforts need to be made to study the behavior of the solution at the endpoints in order to guarantee
the solution derived is within the class of integrable functions. The matrix factorization problem
\beq
G(t)=X^+(t)[X^-(t)]^{-1}=[X^-(t)]^{-1} X^+(t), \quad -k<t<k,
\label{4.7}
\eeq
has the solution 
\beq
X(z)=\GL(z)\left(
\begin{array}{cc}
c(z) \; & \;  s_+(z)\\
s_-(z)\; & \;  c(z)\\
\end{array}\right). 
\label{4.8}
\eeq
Here, we denoted
\beq
c(z)=\cosh[\sqrt{1-z^2}\Gb(z)], \quad s_\pm(z)=(1-z^2)^{\mp 1/2}\sinh[\sqrt{1-z^2}\Gb(z)].
\label{4.9}
\eeq
The function $\sqrt{1-z^2}$ is the single branch of the two-valued function $\Go^2=1-z^2$ 
fixed by the condition $\sqrt{1-z^2}|_{z=0}=1$ in the plane
cut along the straight line joining the branch points $z=-1$ and $z=1$ and passing through the infinite point.
For $-k<t<k$, it is directly verified that
$$
X^+(t)[X^-(t)]^{-1}=\fr{\GL^+(t)}{\GL^-(t)}
$$
\beq
\times\left(
\begin{array}{cc}
\cosh[\sqrt{1-t^2}(\Gb^+(t)-\Gb^-(t))]\; & \;  \fr{1}{\sqrt{1-t^2}}\sinh[\sqrt{1-t^2}(\Gb^+(t)-\Gb^-(t))]\\\
\sqrt{1-t^2}\sinh[\sqrt{1-t^2}(\Gb^+(t)-\Gb^-(t))] \; & \;  \cosh[\sqrt{1-t^2}(\Gb^+(t)-\Gb^-(t))]\\
\end{array}\right). 
\label{4.9'}
\eeq
The matrix $X(z)$ is a solution of the factorization problem if the functions $\GL(z)$ and $\Gb(z)$
solve the following scalar Riemann-Hilbert problems:
\beq
\fr{\GL^+(t)}{\GL^-(t)}=\sqrt{\Gl_1}, \quad -k<t<k,
\label{4.10}
\eeq
and 
\beq
\Gb^+(t)-\Gb^-(t)=\fr{\ln\sqrt{\Gl_1}}{\sqrt{1-t^2}}, \quad -k<t<k, 
\label{4.11}
\eeq
where
\beq
\sqrt{\Gl_1}=i\sqrt{\nu_0}, \quad \ln\sqrt{\Gl_1}=\fr{1}{2}\ln\nu_0+\fr{i\pi}{2}, \quad \nu_0=\fr{3-\nu}{1+\nu}>0.
\label{4.12}
\eeq
The solution of the factorization problem (\ref{4.10}) is defined up  to a rational factor. For our purposes, we choose
it in the form
\beq
\GL(z)=\fr{1}{z-k}\exp\left\{\fr{1}{2\pi i}\int_{-k}^k\fr{\ln\sqrt{\Gl_1}d\tau}{\tau-z}\right\}=(z-k)^{-3/4-i\Gg}(z+k)^{-1/4+i\Gg},
\quad \Gg=\fr{\ln\nu_0}{4\pi}.
\label{4.13}
\eeq
Here, $\GL(z)$ is the branch fixed by the condition $\GL(z)\sim z^{-1}$, $z\to\infty$, in the plane cut along 
the segment $[-k,k]$ passing through the point $z=0$. This branch is positive for $z=t$, $t\in (k,+\infty)$ and negative for
$t\in(-\infty,-k)$. The solution
of the problem (\ref{4.11}) is given by the Cauchy integral
\beq
\Gb(z)=\fr{\ln\sqrt{\Gl_1}}{2\pi i}\int_{-k}^k\fr{d\tau}{\sqrt{1-\tau^2}(\tau-z)}.
\label{4.14} 
\eeq
It has logarithmic singularities at the endpoints, and  the functions $c(z)$ and
$s_\pm(z)$ have power singularities
$$
c(z)\sim C_\pm (z\mp k)^{-1/4+i\Gg}, \quad 
s_+(z)\sim \mp \fr{C_\pm}{\sqrt{1-k^2}} (z\mp k)^{-1/4+i\Gg}, 
$$
\beq
s_-(z)\sim \mp C_\pm\sqrt{1-k^2} (z\mp k)^{-1/4+i\Gg}, \quad z\to\pm k, \quad C_\pm=\const\ne 0.
\label{4.15}
\eeq
In view of the power singularities of the function $\GL(z)$, the matrix $X(z)$ and its inverse
\beq
[X(z)]^{-1}=(z-k)^{3/4+i\Gg}(z+k)^{1/4-i\Gg}\left(
\begin{array}{cc}
c(z) \; & \;  -s_+(z)\\
-s_-(z)\; & \;  c(z)\\
\end{array}\right)
\label{4.16}
\eeq
in neighborhoods of the endpoints behave as
$$
X(z)\sim (z-k)^{-1}Y_+'(z), \quad [X(z)]^{-1}\sim (z-k)^{1/2}Y_+''(z), \quad z\to k,
$$
\beq
X(z)\sim (z+k)^{-1/2}Y_-'(z), \quad [X(z)]^{-1}\sim Y_-''(z), \quad z\to -k,
\label{4.17}
\eeq
where $Y_\pm'(z)$ and $Y_\pm''(z)$ are $2\times 2$ matrices whose elements  are bounded as $z\to\pm k$.

We next study the behavior of the matrices $X(z)$ and $[X(z)]^{-1}$ at the infinite point. 
From (\ref{4.14}) we derive
\beq
\Gb(z)\sim\fr{\Gb_0}{z}, \quad z\to\infty, \quad \Gb_0=\left(2i\Gg-\fr12\right)\sin^{-1} k.
\label{4.19}
\eeq
The branch of the 
function $\sqrt{1-z^2}$ chosen before is discontinuous at the infinite point since the cut passes through
the point $z=\infty$. We have
\beq
\sqrt{1-z^2}\sim -iz\sgn \I z, \quad z\to \infty.
\label{4.18}
\eeq
However, for the functions $c(z)$, $s_+(z)$, and $s_-(z)$, the infinite point is a regular point, a simple zero, and
a simple pole, respectively. 
Combining  this with (\ref{4.19}), (\ref{4.13}), (\ref{4.8}), (\ref{4.9}) and (\ref{4.16}) we find the asymptotics
of the matrix $X(z)$ and its inverse as
\beq
X(z)\sim \fr{1}{z}\left(
\begin{array}{cc}
\cos\Gb_0 \; & \;  z^{-1}\sin\Gb_0\\
-z\sin\Gb_0\; & \;  \cos\Gb_0\\
\end{array}\right),\quad 
[X(z)]^{-1}\sim z\left(
\begin{array}{cc}
\cos\Gb_0 \; & \;  -z^{-1}\sin\Gb_0\\
z\sin\Gb_0\; & \;  \cos\Gb_0\\
\end{array}\right),\quad z\to\infty.
\label{4.20}
\eeq
Substituting the splitting (\ref{4.7}) into the Riemann-Hilbert boundary condition (\ref{4.5}) and replacing
the vector $-i(\nu+1)^{-1}[X^+(t)]^{-1}\BT(t)g(t) $ by $\BGY^+(t)-\BGY^-(t)$
yield
\beq
[X^+(t)]^{-1}\BGF^+(t)-\BGY^+(t)=[X^-(t)]^{-1}\BGF^-(t)-\BGY^-(t), \quad -k<t<k.
\label{4.21}
\eeq
Here, $\BGY^\pm(t)=\BGY(t\pm i0)$ and 
\beq
\BGY(z)=-\fr{1}{2\pi(\nu+1)}\int_{-k}^k\fr{[X^+(\tau)]^{-1}\BT(\tau)g(\tau)d\tau}{\tau-z}.
\label{4.22}
\eeq
Due to (\ref{4.1}) the vector $\BGF(z)$ has a simple zero at the infinite point, and according to the continuity principle, the 
Liouville theorem and the second asymptotic relation in (\ref{4.20}) the vector $[X(z)]^{-1}\BGF(z)-\BGY(z)$
is a polynomial vector of the form
\beq
[X(z)]^{-1}\BGF(z)-\BGY(z)=\left(
\begin{array}{c}
C_0\cos\Gb_0 \\
C_0z\sin\Gb_0 +C_1\\
\end{array}\right),
\label{4.23}
\eeq
where $C_0$ and $C_1$ are arbitrary constants. This asserts 
\beq
\BGF(z)=X(z)\left[\BGY(z)+\left(
\begin{array}{c}
C_0\cos\Gb_0 \\
C_0z\sin\Gb_0 +C_1\\
\end{array}\right)\right], \quad z\in {\Bbb C}\setminus[-k,k].
\label{4.24}
\eeq
Analysis of the asymptotics of the vector $\BGF(z)$ shows that is has an integrable singularity
at the endpoint $z=-k$, while at the second end, $z=k$, it has a non integrable  singularity. It follows 
from (\ref{4.15}) and (\ref{4.24}) that in the vicinity of the point $z=k$
\beq
\BGF(z)=\fr{C_+(z-k)^{-1}}{(2k)^{1/4-i\Gg}}
\left(
\begin{array}{cc}
1 \; & -\fr{1}{\sqrt{1-k^2}} \\
-\sqrt{1-k^2} \; & \; 1\\
\end{array}\right)
\left(
\begin{array}{c}
\GY_1(k)+C_0\cos\Gb_0 \\
\GY_2(k)+C_0k\sin\Gb_0 +C_1\\
\end{array}\right)+\BGF_0(z),
\label{4.25}
\eeq
where the vector $\BGF_0(z)$ may have at most an integrable singularity at the point $z=k$. 
From this relation we infer that the vector $\BGF(z)$ has an integrable singularity at the point $z=k$
if and only if 
\beq
C_1=\sqrt{1-k^2}\GY_1(k)-\GY_2(k)+C_0(\sqrt{1-k^2}\cos\Gb_0-k\sin\Gb_0).
\label{4.26}
\eeq
At the infinite point, both components of the vector $\BGF(z)$ in (\ref{4.24}) have a simple zero and
\beq
\BGF(z)\sim \fr{1}{z}\left(
\begin{array}{c}
C_0  \\
C_1\cos\Gb_0 -\GY_1^\circ\sin\Gb_0\\
\end{array}\right), \quad z\to\infty, 
\label{4.27}
\eeq
where $\GY_1^\circ$ is the first component of the vector
\beq
\BGY^\circ=\fr{1}{2\pi(\nu+1)}\int_{-k}^k[X^+(\tau)]^{-1}\BT(\tau)g(\tau)d\tau.
\label{4.28}
\eeq
Comparing (\ref{4.29}) with (\ref{4.6'}) determines the constant $C_0$
\beq
C_0=-\fr{P_*}{2\pi a}.
\label{4.29}
\eeq
This completes the solution of the vector Riemann-Hilbert problem (\ref{4.5}). Its exact solution is 
given by (\ref{4.24}), (\ref{4.26}) and (\ref{4.29}).

\setcounter{equation}{0}
  
\section{Solution of the singular integral equation}\label{sol} 

The solution of equation (\ref{3.20}) is expressed through the solution of the vector Riemann-Hilbert
problem by the formula $\Gy(at)=\GF_1^+(t)-\GF_1^-(t)$, $-k<t<k$. This section transforms this formula into a form convenient for computations.
According to formulas (\ref{4.8}) and (\ref{4.9}) 
 the limiting values of the matrix $X(z)$ admit the representation
\beq
X^\pm(t)=\left(
\begin{array}{cc}
\Gc^\pm_1(t) \; & \Gc_2^\pm(t) \\
(1-t^2)\Gc_2^\pm(t) \; & \; \Gc_1^\pm(t)\\
\end{array}\right), \quad -k<t<k,
\label{5.1.0}
\eeq
where the functions $\Gc^\pm_j(t)$, $j=1,2$ are to be expressed through the functions $\Gb^\pm(t)$ and $\GL^\pm(t)$.
Utilizing (\ref{4.24}) and (\ref{5.1.0}) we write
$$
\Gy(at)=\GO(t)+[\Gc_2^+(t)-\Gc_2^-(t)]
[\sqrt{1-k^2}\GY_1(k)-\GY_2(k)]
$$
\beq
-\fr{P_*}{2\pi a}\{[\Gc_1^+(t)-\Gc_1^-(t)]\cos\Gb_0+[\Gc_2^+(t)-\Gc_2^-(t)][(t-k)\sin\Gb_0+\sqrt{1-k^2}\cos\Gb_0]\},
\label{5.1.2}
\eeq
where
\beq
\GO(t)=\Gc_1^+(t)\GY_1^+(t)-\Gc_1^-(t)\GY_1^-(t)+\Gc_2^+(t)\GY_2^+(t)-\Gc_2^-(t)\GY_2^-(t).
\label{5.1.3}
\eeq
In order to determine the functions $\Gc_j^\pm(t)$, we evaluate the principal value of the singular integral (\ref{4.14}) that is
\beq
\Gb(t)=\left(\fr14-i\Gg\right)I(t),
\label{5.2}
\eeq
where
\beq
I(t)=\lim_{\Gve\to 0^+}\left(-\int_{-k}^{t-\Gve}\fr{d\tau}{\sqrt{1-\tau^2}(t-\tau)}+\int_{t+\Gve}^k
\fr{d\tau}{\sqrt{1-\tau^2}(\tau-t)}\right).
\label{5.3}
\eeq
Making the substitutions $\Gx=(t-\tau)^{-1}>0$ and $\Gx=(\tau-t)^{-1}>0$ in the first and second integrals
in (\ref{5.3}), respectively, we deduce
\beq
I(t)=-\lim_{\Gve\to 0^+}
\left(\int_{1/(k+t)}^{1/\Gve}\fr{d\Gx}{\sqrt{r_+(\Gx)}}+\int_{1/\Gve}^{1/(k-t)}\fr{d\Gx}{\sqrt{r_-(\Gx)}}
\right),
\label{5.4}
\eeq
where 
\beq
r_\pm(\Gx)=(1-t^2)\Gx^2\pm 2t\Gx-1,\quad 
r_+(\Gx)>0, \quad \fr{1}{k+t}<\Gx<\fr{1}{\Gve}, \quad  r_-(\Gx)>0,
\quad 
\fr{1}{\Gve}<\Gx<\fr{1}{k-t}.
\label{5.4'}
\eeq
Employing the antiderivative
\beq
\int\fr{d\Gx}{\sqrt{r_\pm(\Gx)}}
=\fr{1}{\sqrt{1-t^2}}\ln\left[2\sqrt{(1-t^2)r_\pm(\Gx)}+2(1-t^2)\Gx\pm 2t\right]+\const
\label{5.5}
\eeq
and passing to the limit $\Gve\to 0^+$ in (\ref{5.4}) we eventually obtain
\beq
I(t)=\fr{1}{\sqrt{1-t^2}}\ln\fr{(k-t)R(t)}{k+t},
\label{5.6}
\eeq
where we denoted
$$
R(t)=\fr{R_+(t)}{R_-(t)},
$$
\beq
R_\pm(t)=\sqrt{(1-t^2)[(1-t^2)\pm 2t(k\pm t)-(k\pm t)^2]}+1-t^2\pm t(k\pm t)>0, \quad  -k<t<k.
\label{5.7}
\eeq
We apply now the Sokhotski-Plemelj formulas and discover the following explicit representations
of the limiting values of the function $\Gb(z)$: 
\beq
\sqrt{1-t^2}\Gb^\pm(t)=\left(\fr14-i\Gg\right)\left(\pm \pi i +\ln\fr{(k-t)R(t)}{k+t}\right).
\label{5.8}
\eeq
Utilizing these representations  determines the functions
$$
\cosh[\sqrt{1-t^2}\Gb^\pm(t)]=\fr{e_0^{\pm 1}}{2}
\left[\fr{(k-t)R(t)}{k+t}\right]^{1/4-i\Gg}+
\fr{e_0^{\mp 1}}{2}\left[\fr{(k-t)R(t)}{k+t}\right]^{-1/4+i\Gg},
$$
\beq
\sinh[\sqrt{1-t^2}\Gb^\pm(t)]=\fr{e_0^{\pm 1}}{2}\left[\fr{(k-t)R(t)}{k+t}\right]^{1/4-i\Gg}-
\fr{e_0^{\mp 1}}{2}\left[\fr{(k-t)R(t)}{k+t}\right]^{-1/4+i\Gg},
\label{5.9}
\eeq
Here, $e_0=e^{\pi i/4+\pi \Gg}.$
The limiting values of the function $\GL(z)$ are discontinuous through the cut $[-k,k]$. According to the 
choice of the branch of $\GL(z)$, $\arg(z-k)=\pi$ and $\arg(z+k)=\pi\mp\pi$ as $z=t\pm i0$, $-k<t<k$. 
We have
\beq
\GL^\pm(t)=-e_0^{\pm 1}(k-t)^{-3/4-i\Gg}(k+t)^{-1/4+i\Gg}, \quad -k<t<k.
\label{5.10}
\eeq
Combining these results we find the limiting values $\Gc_j^\pm(t)$, $-k<t<k$, of the functions  $\Gc_j(z)$. 
They are  given by
\beq
\Gc_j^\pm(t)=-\fr{1}{2(\sqrt{1-t^2})^{j-1}}\left[e_0^{\pm 2}(k-t)^{-1/2-2i\Gg}(k+t)^{-1/2+2i\Gg}R^{1/4-i\Gg}(t)
-\fr{(-1)^j}{k-t}R^{-1/4+i\Gg}(t)\right].
\label{5.11}
\eeq
Now we wish to specify the limiting values of the vector $\BGY(z)$. Computing first 
\beq
[X^+(t)]^{-1}\BT(t)=-\fr{1}{e_0^2}(k-t)^{1/2+2i\Gg}(k+t)^{1/2-2i\Gg}R^{-1/4+i\Gg}\BT(t)
\label{5.12}
\eeq
and employing next the Sokhotski-Plemelj formulas we obtain the components $\GY_1^\pm(t)$
and $\GY_2^\pm(t)$ of the vectors $\BGY^\pm(t)$  as
\beq
\GY_j^\pm(t)=\pm\fr{ig(t)}{2(\nu+1)e_0^2}(\sqrt{1-t^2})^{j-1}(k-t)^{1/2+2i\Gg}(k+t)^{1/2-2i\Gg}R^{-1/4+i\Gg}(t)+\GY_j(t),
\label{5.13}
\eeq
where $\GY_j(t)$ are the principal values of the integrals
\beq
\GY_j(t)=\fr{1}{2\pi(\nu+1)e_0^2}\int_{-k}^k
(k-\tau)^{1/2+2i\Gg}(k+\tau)^{1/2-2i\Gg}R^{-1/4+i\Gg}(\tau)\fr{(\sqrt{1-\tau^2})^{j-1}g(\tau)d\tau}{\tau-t}.
\label{5.14}
\eeq
It is seen that the functions $\GY_j^\pm(t)$ are bounded at the endpoints and have finite limits. At the point $k$,
\beq
\GY_j(k)= -\fr{1}{2\pi(\nu+1)e_0^2}\int_{-k}^k
(k-\tau)^{-1/2+2i\Gg}(k+\tau)^{1/2-2i\Gg}R^{-1/4+i\Gg}(\tau)(\sqrt{1-\tau^2})^{j-1}g(\tau)d\tau.
\label{5.15}
\eeq
We now ready to evaluate  $\GO(t)$ and $\Gc_j^+(t)-\Gc_j^-(t)$. 
Substituting
(\ref{5.11}) and (\ref{5.13}) into (\ref{5.1.3}) yields, after some simplifications,
$$
\GO(t)=-\fr{i(1-\nu)g(t)}{\nu_1}+\Gc(t)\left[\GY_1(t)+\fr{\GY_2(t)}{\sqrt{1-t^2}}\right],
$$
\beq
\Gc_j^+(t)-\Gc_j^-(t)=\fr{\Gc(t)}{(\sqrt{1-t^2})^{j-1}}, \quad j=1,2.
\label{5.16}
\eeq
Here,
\beq
\Gc(t)=-\fr{2i}{(1+\nu)\sqrt{\nu_0}}(k-t)^{-1/2-2i\Gg}(k+t)^{-1/2+2i\Gg}R^{1/4-i\Gg}(t).
\label{5.17}
\eeq
Remembering the formula for the function $\psi(at)$, we 
can write down the final  representation for the solution of the singular integral equation (\ref{3.20})
$$
\Gy(at)=-\fr{i(1-\nu)g(t)}{\nu_1}
+\Gc(t)\left[\GY_1(t)+\sqrt{\fr{1-k^2}{1-t^2}}\GY_1(k)
\right.
$$
\beq
\left.
+\fr{\GY_2(t)-\GY_2(k)}{\sqrt{1-t^2}}-\fr{P_*}{2\pi a}\left(\cos\Gb_0+\fr{(t-k)\sin\Gb_0+\sqrt{1-k^2}\cos\Gb_0}
{\sqrt{1-t^2}}\right)\right].
\label{5.18}
\eeq
Notice that this formula is a closed-form solution  and possesses  only two integrals (\ref{5.14}).

\setcounter{equation}{0}
  
\section{Analysis of the solution at the singular points. Numerical results}\label{ana}

\subsection{Stress intensity factors}

In Section \ref{SIE}, we saw that  the jump
of the tangential derivative of the displacement vector could readily be expressed through the solution of equation (\ref{3.20}) by (\ref{3.19'}). Due to formulas (\ref{5.17}) and (\ref{5.18}) this function has the square root singularities
at the crack tips
\beq
\Gf(x)\sim\Gf_0^\pm(a\mp x)^{-1/2}, \quad x\to \pm a^\mp,
\label{6.3}
\eeq
 where $a^\mp=a\mp 0$ and
 \beq
 \Gf_0^\pm=\mp\fr{1}{\pi}\sqrt{\fr{2}{a}}\left[\int_{-b}^b\left(\fr{a+\Gx}{a-\Gx}\right)^{\pm 1/2}\Gy(\Gx)d\Gx+
 2\int_{-a}^a\left(\fr{a+\Gx}{a-\Gx}\right)^{\pm 1/2}\Gs^\circ(\Gx)d\Gx\pm\fr{(1-\nu)P_*}{2}\right].
 \label{6.4}
 \eeq
At the inclusion endpoints, $x=b$ and $x=-b$, the internal singular points, the function $\Gf(x)$
not only has square root singularities but also oscillates
\beq
\Gf(x)=(x-b)^{-1/2-2i\Gg}(x+b)^{-1/2+2i\Gg}\GU(x),
\label{6.4'}
\eeq
where $\GU(x)$ is bounded and in general discontinuos at the points $x=\pm b$. It has definite nonzero limits
as $x\to b^\pm$ and as $x\to -b^\mp$.

Now we wish to analyze the concentration of stresses in the vicinities of the crack and
the inclusion tips.  Assume $y=0$ and $x$ is outside the segment $[-b,b]$. Applying (\ref{3.11})
we have
\beq
\Gs_{12}(x,0)+i\Gs_{22}(x,0)=
\fr{i(1-\nu)}{4\pi}\int_{-b}^b\fr{\psi(\Gx)d\Gx}{\Gx-x}+\fr{1}{4\pi}\int_{-a}^a
\fr{\Gf(\Gx)d\Gx}{\Gx-x},
\quad |x|>b.
\label{6.2}
\eeq
 Analysis of the Cauchy integral with the density $\Gf(\Gx)$ in (\ref{6.2}) yields
 \beq
 \Gs_{12}(x,0)+i\Gs_{22}(x,0)\sim\mp\fr{\Gf_0^\pm}{4\sqrt{\pm x-a}}, \quad x\to \pm a^\pm.
 \label{6.5}
 \eeq
 On the other hand, on employing the conventional notations of the stress intensity factors (SIFs)
 we have
 \beq
 \Gs_{12}(x,0)+i\Gs_{22}(x,0)\sim\fr{K^\pm_{II}+iK_I^\pm}{\sqrt{2\pi(\pm x-a)}}, \quad  x\to \pm a^\pm.
\label{6.6}
\eeq
Combining  (\ref{6.5}) and (\ref{6.6}) we express the SIFs through two integrals 
\beq
K^\pm_{II}+iK_I^\pm=\fr{1}{2\sqrt{\pi a}}
\left[\int_{-b}^b\left(\fr{a+\Gx}{a-\Gx}\right)^{\pm 1/2}\Gy(\Gx)d\Gx+
 2\int_{-a}^a\left(\fr{a+\Gx}{a-\Gx}\right)^{\pm 1/2}\Gs^\circ(\Gx)d\Gx\pm\fr{(1-\nu)P_*}{2}\right].
 \label{6.7}
 \eeq
Next we  determine  the contact stresses (the traction vector components) 
\beq
 \Gs_{12}(x,0^-)+i\Gs_{22}(x,0^-)=\Gs_+(x)-\psi(x)
 \label{6.8}
 \eeq
acting in the contact zone $-b<x<b$.  After substituting (\ref{5.18}) and (\ref{3.21}) into (\ref{6.8}) and since $\Gs_+(x)=-\Gs^\circ(x)$,
we have 
 $$
 \Gs_{12}(at,0^-)+i\Gs_{22}(at,0^-)=
- \fr{4\Gs^\circ(at)}{\nu_1}+\fr{1-\nu}{\nu_1}\left[
E h'(at)+iEw^\circ(at)
 \right.
 $$$$
 \left.
 -\fr{2i}{\pi\sqrt{1-t^2}}\
 \int_{-1}^1
\fr{\sqrt{1-\tau^2}}{\tau-t}\Gs^\circ(a\tau)d\tau
+\fr{i(1-\nu)P_*}{2\pi a\sqrt{1-t^2}}\right]
-\Gc(t)\left[\GY_1(t)+\sqrt{\fr{1-k^2}{1-t^2}}\GY_1(k)
\right.
$$
\beq
\left.
+\fr{\GY_2(t)-\GY_2(k)}{\sqrt{1-t^2}}-\fr{P_*}{2\pi a}\left(\cos\Gb_0+\fr{(t-k)\sin\Gb_0+\sqrt{1-k^2}\cos\Gb_0}
{\sqrt{1-t^2}}\right)\right], \quad -k<t<k.
\label{6.9}
\eeq
By examining this formula we see that due to the existence of the singularities of the function $\Gc(t)$ in the vicinity of the inclusion tips, $t\to\pm k^\mp$ ($x\to \pm b^\mp)$, the traction admits the representation
\beq
 \Gs_{12}(at,0^-)+i\Gs_{22}(at,0^-)=(k-t)^{-1/2-2i\Gg}(k+t)^{-1/2+2i\Gg}\Xi(t),
 \label{6.10}
 \eeq
 where $\Xi(t)$ is a bounded function having definite limits as $t\to \pm k^\mp.$ 
 
 \begin{figure}[t]
\centerline{
\scalebox{0.6}{\includegraphics{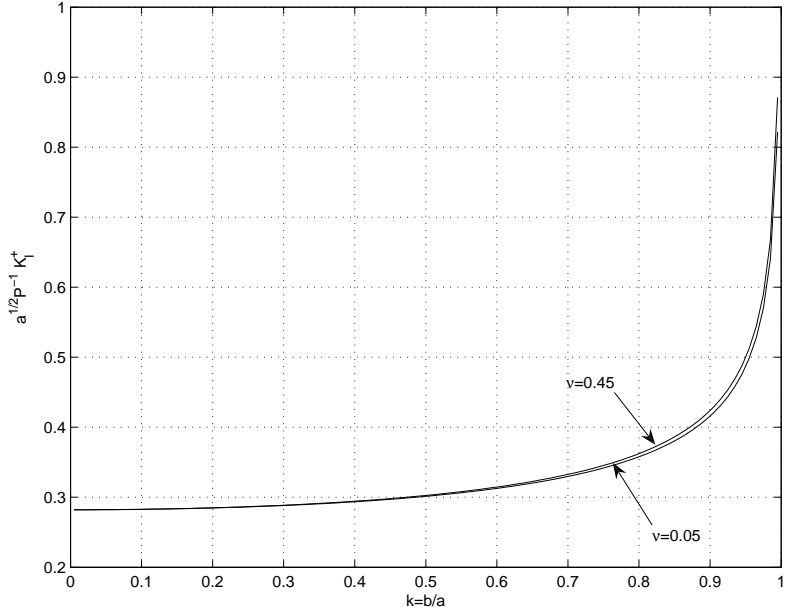}}}
\caption{The SIF $K_I^+$ at the crack tip $x=a$ vs $k=b/a$ for $\nu=0.1$ and $\nu=0.5.$}
\label{fig2}
\end{figure} 

\begin{figure}[t]
\centerline{
\scalebox{0.6}{\includegraphics{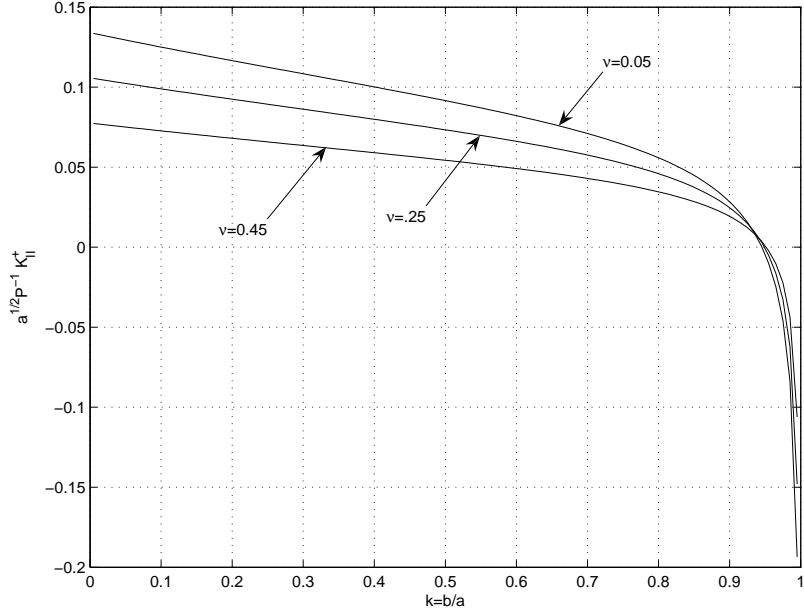}}}
\caption{The SIF $K_{II}^+$ at the crack tip $x=a$ vs $k=b/a$ for $\nu=0.1$, $\nu=0.3$, and $\nu=0.5.$}
\label{fig3}
\end{figure}

\begin{figure}[t]
\centerline{
\scalebox{0.6}{\includegraphics{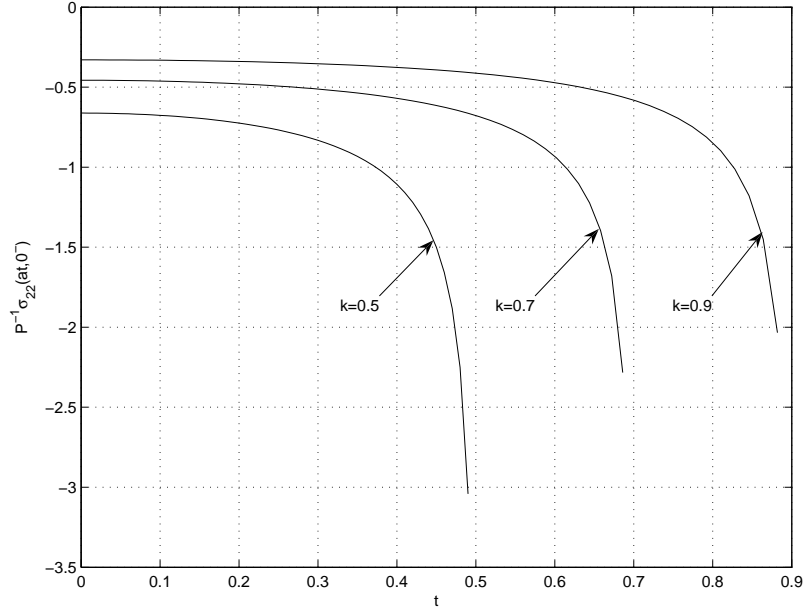}}}
\caption{The normal traction component  $P^{-1}\Gs_{22}(at,0^-)$  vs $t$ for  $k=0.5$, $k=0.7$,  and $k=0.9$ 
when  $\nu=0.3.$}
\label{fig4}
\end{figure}

\begin{figure}[t]
\centerline{
\scalebox{0.6}{\includegraphics{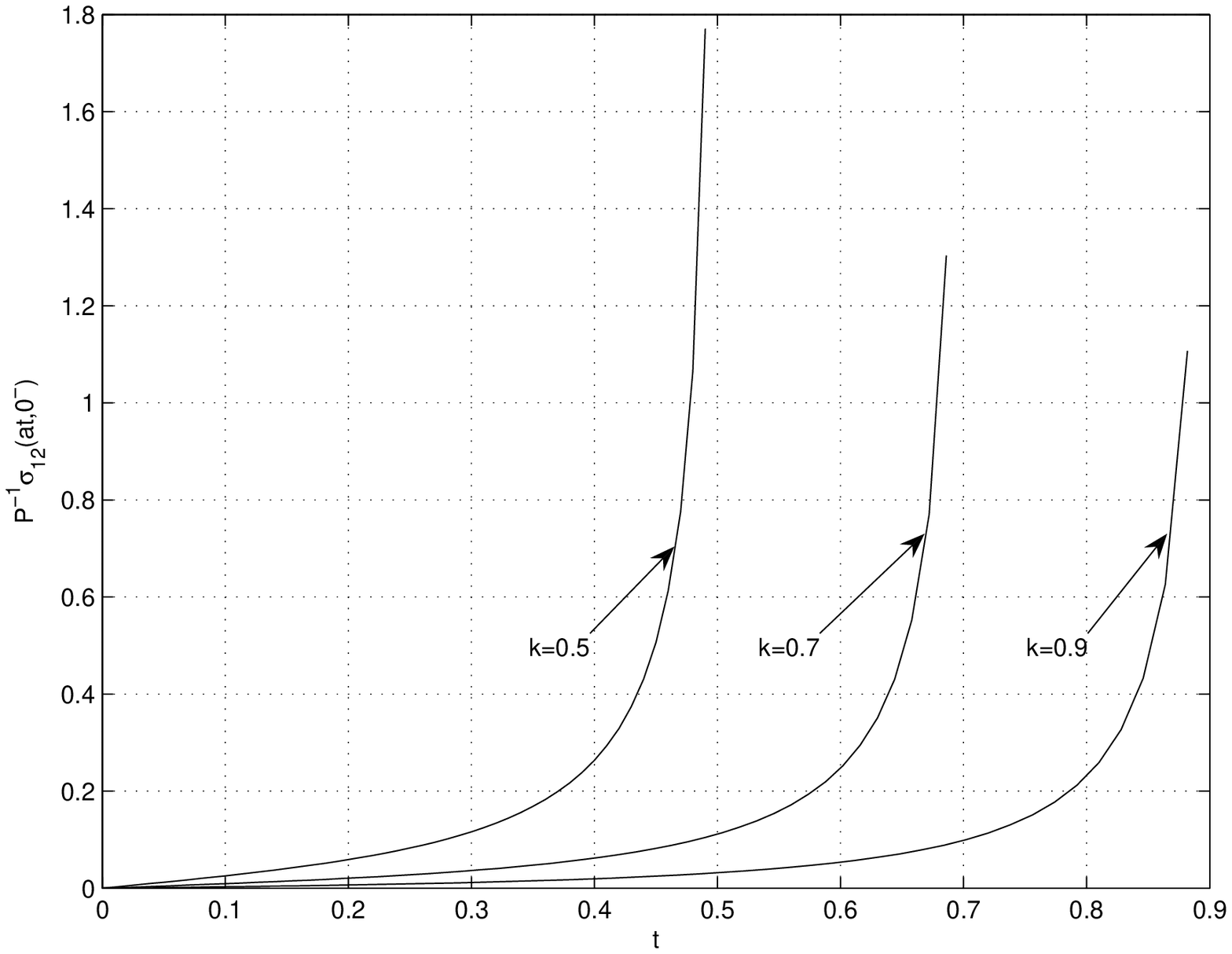}}}
\caption{
The tangential traction component  $P^{-1}\Gs_{12}(at,0^-)$  vs $t$  for  $k=0.5$, $k=0.7$,  and $k=0.9$ 
when  $\nu=0.3.$}
\label{fig5}
\end{figure}

\subsection{Computational formulas for the SIFs}

For a numerical example we take  $w^\circ(x)=0$, $h'(x)=0$, $-b<x<b$, and  $\Gs^\circ(x)=0$, $-a<x<a$. Then
the function $g(t)$ has a simple form,
\beq
g(t)=\fr{(1-\nu)P}{2\pi a \sqrt{1-t^2}}.
\label{6.11}
\eeq
To evaluate the function $\psi(at)$ given by (\ref{5.18}) we need to compute the functions $\Psi_j(t)$.
It is convenient to represent them as
\beq
\Psi_j(t)=P_0\int_{-k}^k\sqrt{k^2-\tau^2}g_j(\tau)\fr{d\tau}{\tau-t},\quad -k<t<k, \quad j=1,2,
\label{6.12}
\eeq
where
\beq
g_1(\tau)=\left(\fr{k-\tau}{k+\tau}\right)^{2i\Gg}\fr{R^{-1/4+i\Gg}(\tau)}{\sqrt{1-\tau^2}}, \quad
g_2(\tau)=\sqrt{1-\tau^2}g_1(\tau),\quad
P_0=\fr{(1-\nu)P}{4\pi^2 a(1+\nu)e_0^2}.
\label{6.13}
\eeq
The functions $g_j(\tau)$ can be expanded in terms of the Chebyshev polynomials of the second kind as
\beq
g_j(\tau)=\sum_{m=0}^\infty c_{jm}U_m\left(\fr{\tau}{k}\right),\quad -k<t<k,
\label{6.14}
\eeq
with the coefficients $c_{jm}$  defined by
\beq
c_{jm}=\fr{2}{\pi k^2}\int_{-k}^k g_j(\tau)\sqrt{k^2-\tau^2}U_m\left(\fr{\tau}{k}\right)d\tau.
\label{6.15}
\eeq
On employing next the spectral relations for the Chebyshev polynomials of the second and first kind
\beq
\int_{-k}^k \sqrt{k^2-\tau^2}U_m\left(\fr{\tau}{k}\right)\fr{d\tau}{\tau-t}=-\pi k T_{m+1}\left(\fr{t}{k}\right),
\quad -k<t<k, \quad  m=0,1,\ldots,
\label{6.16}
\eeq
we obtain the series representation of the functions $\Psi_j(t)$
\beq
\Psi_j(t)=-\pi k P_0\sum_{m=0}^\infty c_{jm} T_{m+1}\left(\fr{t}{k}\right),\quad -k<t<k.
\label{6.17}
\eeq
Since the functions $g_j(\tau)$ oscillate at the points $\tau=\pm k$, instead of using formula (\ref{6.17})
we evaluate $\Psi_j(k)$ directly
\beq
\Psi_j(k)=-P_0\int_{-k} ^k 
\left(\fr{k+\tau}{k-\tau}\right)^{1/2} 
g_j(\tau)d\tau.
\label{6.18}
\eeq
To write down the final formula used for the SIFs, we denote
$$
\Gc_\pm(t)=-\fr{-2i a}{(1+\nu)\sqrt{\nu_0}g_2(t)}\left(\fr{1+t}{1-t}\right)^{\pm 1/2} 
\left[\GY_1(t)+\sqrt{\fr{1-k^2}{1-t^2}}\GY_1(k)
\right.
$$
\beq
\left.
+\fr{\GY_2(t)-\GY_2(k)}{\sqrt{1-t^2}}-\fr{P}{2\pi a}\left(\cos\Gb_0+\fr{(t-k)\sin\Gb_0+\sqrt{1-k^2}\cos\Gb_0}
{\sqrt{1-t^2}}\right)\right].
\label{6.19}
\eeq
In terms of these functions, after a rearrangement, formula (\ref{6.7}) reads
\beq
K^\pm_{II}+iK_I^\pm=\fr{1}{2\sqrt{\pi a}}
\left[\pm\fr{(1-\nu)P}{2}-\fr{i(1-\nu)^2P}{2\pi \nu_1}\ln \fr{1+k}{1-k}+
\int_{-k}^k\fr{\Gc_\pm(\tau)d\tau}{\sqrt{k^2-\tau^2}}\right].
 \label{6.20}
 \eeq
The integrals (\ref{6.15}), (\ref{6.18}), and (\ref{6.20}) are computed by employing the corresponding
Gauss quadrature formulas
$$
c_{jm}\approx \fr{2}{N+1}\sum_{l=1}^N\sin\fr{l\pi}{N+1}\sin\fr{l\pi(m+1)}{N+1}g_j\left(k\cos\fr{l\pi}{N+1}\right),
$$$$
\Psi_j(k)\approx -\fr{4\pi kP_0}{2M+1}\sum_{l=1}^M x_lg_j(-k+2k x_l), \quad x_l=\cos^2\fr{(2l-1)\pi}{2(2M+1)},
$$
\beq
\int_{-k}^k \fr{\Gc_\pm(\tau)d\tau}{\sqrt{k^2-\tau^2}}\approx \fr{\pi}{L}\sum_{l=1}^L\Gc_\pm\left(k\cos\fr{(2l-1)\pi}{2 L}\right),
\label{6.21}
\eeq
where $N$, $M$, and $L$ are the orders of the Gauss integration formulas.

 Figures 2 and 3  represent the dimensionless SIFs $\hat K_I^+=\sqrt{a}P^{-1} K_I$
and $\hat K_{II}^+=\sqrt{a}P^{-1} K_{II}$ at the crack tip $x=a$ for some values of the Poisson ratio $\nu$.
At the left crack tip, $x=-a$, for the symmetric case (\ref{6.11}), $K^-_I=K^+_I$, $K^-_{II}=-K^+_{II}$.
Referring to these figures, we  observe the effect of the ratio $k=b/a$ on the SIFs. As the inclusion ends
approach the crack tips, the Mode-I SIFs $K_I^\pm$ grow to $+\infty$, while the Mode-II SIF $K_{II}^+$ decreases to zero for $k$ close to $0.95$ and then, as $k\to 1$, $K_{II}^+\to-\infty$.
Figures 4 and 5 show sample curves of the traction components in the right-hand contact zone
$0<x<b$, $y=0^-$. Due to the symmetry, $\Gs_{12}(x,0^-)=-\Gs_{12}(-x,0^-)$ and $\Gs_{12}(x,0^-)=\Gs_{22}(-x,0^-)$.
It is seen that the absolute values of both contact stresses attain minima at the point $x=0$ and 
grow as $x\to b^-$ ($t\to k^-$).

\subsection{Limiting case $k\to 0$}

To test the numerical results obtained, we derive explicit formulas for the SIFs in the limiting case 
$k\to 0$ by passing to the limit in (\ref{6.7}). In the case (\ref{6.11}) formula (\ref{6.7}) becomes
\beq
K^\pm_{II}+iK_I^\pm=\fr{1}{2\sqrt{\pi a}}
\left[a\int_{-1}^1\left(\fr{1+k\Gx}{1-k\Gx}\right)^{\pm 1/2}k\Gy(b\Gx)d\Gx\pm
 \fr{(1-\nu)P}{2}\right],
 \label{6.22}
 \eeq
where due to (\ref{5.18}) the function $k\Gy(b\Gx)$ has the form
$$
k\Gy(b\Gx)=-\fr{i(1-\nu)kg(k\Gx)}{\nu_1}
+\Gc^\circ(\Gx)\left[\GY_1(k\Gx)+\sqrt{\fr{1-k^2}{1-k^2\Gx^2}}\GY_1(k)
\right.
$$
\beq
\left.
+\fr{\GY_2(k\Gx)-\GY_2(k)}{\sqrt{1-k^2\Gx^2}}-\fr{P}{2\pi a}\left(\cos\Gb_0+\fr{k(\Gx-1)\sin\Gb_0+\sqrt{1-k^2}\cos\Gb_0}
{\sqrt{1-k^2\Gx^2}}\right)\right],
\label{6.23}
\eeq
and 
\beq
\Gc^\circ(\Gx)=-\fr{2i}{(1+\nu)\sqrt{\nu_0}}(1-\Gx)^{-1/2-2i\Gg}(1+\Gx)^{-1/2+2i\Gg}R^{1/4-i\Gg}(k\Gx).
\label{6.24}
\eeq
It is directly verified that $\Gb_0\to 1$ and $\Psi_j(k)\to 0$ as $k\to 0^+$ and also
that the functions $R(k\Gx)$ and $\Psi_j(k\Gx)$ uniformly with respect to $\Gx\in[-1,1]$
tend to $1$ and $0$, respectively. Passing to the limit $k\to 0^+$ in (\ref{6.22}) we deduce
\beq
\lim_{k\to 0^+}(K^\pm_{II}+iK_I^\pm)=\fr{iP}{\pi\sqrt{\pi a\nu_0}(1+\nu)}\int_{-1}^1
(1-\Gx)^{-1/2-2i\Gg} (1+\Gx)^{-1/2+2i\Gg}d\Gx \pm
 \fr{(1-\nu)P}{4\sqrt{\pi a}},
\label{6.25}
\eeq
This ultimately yields
that the normalized SIFs $K_I^\pm$ in the limiting case $k\to 0$ are independent of the Poisson ratio and  given by
\beq
\sqrt{a}P^{-1}K_I^\pm=\fr{1}{2\sqrt{\pi}}= 0.2820948\ldots.
\label{6.26}
\eeq
The SIFs $K_{II}^\pm$ as $k\to 0$ become
\beq
\sqrt{a}P^{-1}K_{II}^\pm=\pm\fr{1-\nu}{4\sqrt{\pi}}.
\label{6.27}
\eeq
Formulas (\ref{6.26}) and (\ref{6.27}) coincide with the formulas for the SIFs 
derived in  {\bf(\ref{cher2})} for a crack $\{-a<x<a, y=0^\pm\}$
whose upper side is free of traction and the lower side is subjected 
to the normal concentrated load $P$ applied at the center, $\Gs_{12}(x,0^-)=0$, $\Gs_{22}(x,0^-)=-P\Gd(x)$.

Notice that formulas (\ref{6.26}) and (\ref{6.27}) are consistent with the values of the SIF $K_I^+$ 
and $K_{II}^+$ for small values of the parameter $k$ obtained
from formula (\ref{6.7}) in the case (\ref{6.11}): if $k=0.005$, then  $\sqrt{a}P^{-1}K_I^+=0.2820933$ for $\nu=0.05$
and  $\sqrt{a}P^{-1}K_I^+=0.2821035$ for $\nu=0.45$.
For  $\nu=0.05$, $\nu=0.25$, and $\nu=0.45$  formula (\ref{6.27}) gives the following
values of the dimensionless SIF $\sqrt{a}P^{-1}K_{II}^+$: 0.1339950, 0.1057855, and
0.07757607, respectively. On the other hand, the corresponding values
of this factor for $k=0.005$ computed from formula  (\ref{6.7}) in the case (\ref{6.11}) are
0.1335321, 0.1054322, and 0.07732310; they are in good agreement with the limiting values for $k=0$.

\begin{figure}[t]
\centerline{
\scalebox{0.6}{\includegraphics{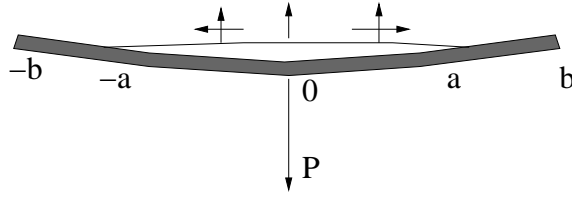}}
}
\caption{The geometry of Model 2 }
\label{fig6}
\end{figure}

\setcounter{equation}{0}
  
\section{Cases $a<b$ and $a=b$}\label{cases}

\subsection{An inclusion whose upper surface is partly separated from the matrix: $a<b$} 

Suppose now that the crack length is less than that of the inclusion that is $a<b$ (Fig. 6).
In the notations (\ref{3.10}), the boundary conditions 
(\ref{2.1}) should be replaced by
$$
w_-(x)=ih'(x)-w^\circ(x), \quad -b<x<b,
$$$$
w_+(x)=ih'(x)-w^\circ(x), \quad a<|x|<b,
$$
\beq
\Gs_+(x)=-\Gs^\circ(x), \quad -a<x<a.
\label{7.2}
\eeq
We will show that this problem reduces to a single singular integral equation 
with respect to the function $\Gf(x)$ in the segment $(-a,a)$, and the structure
of the equation is the same as that of equation (\ref{3.20}).
Rewrite first the representations (\ref{3.12}) as
$$
\fr{1}{\pi}\int_{-b}^b\fr{\psi(\Gx)d\Gx}{\Gx-x}=-\fr{2}{\nu_1}[\Gf(x)+2Ew_-(x)]-
\fr{i(1-\nu)}{\pi \nu_1}
\int_{-a}^a\fr{\Gf(\Gx)d\Gx}{\Gx-x}, \quad -b<x<b,
$$
\beq
2\Gs^\circ(x)=-\psi(x)+\fr{i(1-\nu)}{\nu_1}[\Gf(x)+2Ew_-(x)]-\fr{2}{\pi\nu_1}\int_{-a}^a\fr{\Gf(\Gx)d\Gx}{\Gx-x}, \quad
-a<x<a,
\label{7.3}
\eeq
where $\nu_1=(3-\nu)(1+\nu)$ is the parameter introduced in (\ref{3.9'}).
Inverting the Cauchy integral in the left-hand side of the first equation in (\ref{7.3}) and
applying the Poincar\'e-Bertrand formula and the inclusion  equilibrium condition   (\ref{3.22})
similarly to Section \ref{SIE} we express the function $\psi(x)$
through $\Gf(x)$. We have
\beq
\psi(x)=\fr{2}{\pi\nu_1\sqrt{b^2-x^2}}
\int_{-b}^b\fr{\sqrt{b^2-\Gx^2}}{\Gx-x}[\Gf(\Gx)+2Ew_-(\Gx)]d\Gx-\fr{i(1-\nu)}{\nu_1}\Gf(x)+\fr{iP_*}{\pi\sqrt{b^2-x^2}}.
\label{7.5}
\eeq
Upon substituting (\ref{7.5}) into the second equation in (\ref{7.3}) we eventually
derive the governing integral equation for the function $\Gf(x)$. It reads
\beq
\fr{1}{\pi}\int_{-1/k}^{1/k}
\left(1+\sqrt{\fr{1-\tau^2}{1-t^2}}
\right)\fr{\Gf(b\tau)d\tau}{\tau-t}
-i(1-\nu)\Gf(bt)
=g_1(t), \quad -\fr1k<t<\fr1k,
\label{7.6}
\eeq
where $1/k=a/b<1$,
$$
g_1(t)=-\nu_1\Gs^\circ(bt)-E(1-\nu)[h'(bt)+iw^\circ(bt)]$$
\beq
-\fr{2E}{\pi\sqrt{1-t^2}}\int_{-1}^1
\fr{\sqrt{1-\tau^2}}{\tau-t}[ih'(b\tau)-w^\circ(b\tau)]d\tau
-\fr{iP_*\nu_1}{2\pi b\sqrt{1-t^2}}.
\label{7.7}
\eeq
The function $\Gf(at)$  has to satisfy the additional condition (\ref{3.15}).
It is seen that the governing integral equation (\ref{7.6}) coincides with equation (\ref{3.20}) 
 if we replace $\psi(at)$ by $\Gf(bt)$, $k$ by $1/k$,  and $g(t)$ by $g_1(t)$.
 Instead of the additional condition (\ref{3.22}) for the function $\psi$  we now have
 the condition (\ref{3.15}) for the function $\Gf$.

\subsection{An inclusion whose upper surface is completely separated from the matrix: passage to the limit $k\to 1$} 

Consider now the case when the whole upper side of the inclusion is debonded from the matrix
that is when there is a crack $-a<x<a$, $y=0^+$ between the inclusion $-a<x<a$, $y=0^-$ and the matrix $\Bbb R^2\setminus
\{-a<x<a, y=0\}$.
The governing equation in this case can be obtained from (\ref{3.20}) 
\beq
\fr{1}{\pi}\int_{-1}^1\left(1+\sqrt{\fr{1-\tau^2}{1-t^2}}
\right)\fr{\Gy(a\tau)d\tau}{\tau-t}
-i(1-\nu)\Gy(at)=g(t), \quad -1<t<1,
\label{7.8}
\eeq
or from (\ref{7.6}) by putting $k=1$
\beq
\fr{1}{\pi}\int_{-1}^{1}
\left(1+\sqrt{\fr{1-\tau^2}{1-t^2}}
\right)\fr{\Gf(a\tau)d\tau}{\tau-t}
-i(1-\nu)\Gf(at)
=g_1(t), \quad -1<t<1.
\label{7.8'}
\eeq
The former equation, when solved, provides the jump of the traction vector, while
the latter yields the function $\Gf(x)$, the jump of the tangential derivative
of the displacement vector. There are other ways to solve the problem in this case. Before we implement passage
to the limit $k\to 1$ in the solution derived in Sections \ref{RH} and \ref{sol} we briefly describe  a
method for the system
(\ref{3.12}) based on its decoupling.
 Since the traction vector is prescribed in the whole segment $-a<x<a$, 
$\Gs_+(x)=-\Gs^\circ(x)$, and the displacement vector is known on the whole lower side of the inclusion,
$w_-(x)=ih'(x)-w^\circ(x)$,  $-a<x<a$,  the relations (\ref{3.12}) allow us to deduce the following governing system: 
\beq
\BGo(x)+\fr{A}{2\pi}\int_{-a}^a
\fr{\BGo(\Gx)d\Gx}{\Gx-x}=\Bg(x), \quad -a<x<a,
\label{7.9}
\eeq
where
\beq
A=\left(
\begin{array}{cc}
i(1-\nu) \; & \;  1\\
\nu_1 \; & \;  i(1-\nu)\\
\end{array}\right), \quad
\BGo(x)= \left(
\begin{array}{c}
\psi(x) \\
\Gf(x)\\
\end{array}\right), \quad \Bg(x)= -2\left(
\begin{array}{c}
\Gs^\circ(x)(x) \\
E[ih'(x)-w^\circ(x)]\\
\end{array}\right).
\label{7.10}
\eeq
The matrix $A$ has constant entries, and it is diagonalizable, $A=BDB^{-1}$, where $D$ is a diagonal matrix, and 
$B$ is a nonsingular matrix of transformation.
Therefore, the system
(\ref{7.9}) reduces to two scalar singular integral equations which admit a closed-form solution by reducing them to the associated 
scalar Riemann-Hilbert problems in the segment $-1<t<1$. This method was applied for example in  {\bf(\ref{pop})} to solve the problem
in the case $a=b$ for an inclusion in the interface of a composite plane.
Here, we give an alternative solution by solving the integral equation (\ref{7.8}) by 
passing to the limit $k\to 1$ in the representation formulas for the solution in the case $k\in(0,1)$.

Following the scheme of Section \ref{RH} we reduce the singular integral equation ({\ref{7.8}) to the vector Riemann-Hilbert problem
(\ref{4.5}), (\ref{4.6}) on the interval $(-1,1)$ and factorize the matrix $G(t)$ by formulas (\ref{4.7}) to (\ref{4.9}) with the functions
$\GL(z)$ and $\Gb(z)$ being the solutions of the scalar Riemann problems (\ref{4.10}) and (\ref{4.11}), respectively.
The solution of the former problem is the one given by (\ref{4.13}) with $k=1$,
\beq
\GL(z)=(z-1)^{-3/4-i\Gg}(z+1)^{-1/4+i\Gg},
\quad \Gg=\fr{\ln\nu_0}{4\pi},
\label{7.11}
\eeq
and  the single branch is fixed by the condition $\GL(z)\sim z^{-1}$, $z\to\infty$, in the plane cut along 
the segment $[-1,1]$ passing through the point $z=0$.  
In particular, when $z=t\pm i0$ and $-1<t<1$, we have
\beq
\GL^\pm(t)=-e_0^{\pm 1}(1-t)^{-3/4-i\Gg}(1+t)^{-1/4+i\Gg},\quad -1<t<1,
\label{7.12}
\eeq
$\arg(1\pm t)=0$, $t\in(-1,1)$, and $e_0=e^{\pi i/4+\pi \Gg}$ is the parameter introduced in Section \ref{sol}.

Let now $k\to 1^-$ in (\ref{5.8}). Since $R(t)\to (1+t)/(1-t)$ as $k\to 1^-$ and $-1<t<1$, we deduce
\beq
\sqrt{1-t^2}\Gb^\pm(t)=\pm \pi i \left(\fr14-i\Gg\right),
\label{7.13}
\eeq
and the Wiener-Hopf matrix-factors in (\ref{4.7}) become 
\beq
X^\pm(t)=\GL^\pm(t)\left(
\begin{array}{cc}
e_+ \; & \;  \pm e_-(1-t^2)^{-1/2}\\
\pm e_-(1-t^2)^{1/2}\; & \;  e_+\\
\end{array}\right), \quad -1<t<1, \quad e_\pm=\fr12(e_0\pm e_0^{-1}).
\label{7.14}
\eeq
We assert next that as $k\to 1^-$ in (\ref{4.19}),  $\Gb_0\to \pi(-\fr14+i\Gg)$, $\cos\Gb_0\to e_+$,
$\sin\Gb_0\to ie_-$, and therefore the matrices $[X(z)]^{\pm 1}$ 
have the following asymptotics  at infinity:
\beq
[X(z)]^{\pm 1}\sim z^{\mp 1} \left(
\begin{array}{cc}
e_+ \; & \;  \pm z^{-1}ie_-\\
\mp zie_-\; & \;  e_+\\
\end{array}\right),\quad z\to\infty.
\label{7.15}
\eeq
Since the behavior of the matrix $X(z)$ at the infinite point is the same as before, the solution of the vector
Riemann-Hilbert problem is given by (\ref{4.24}) that reads in the limiting case $k\to 1^-$
\beq
\BGF(z)=X(z)\left[\BGY(z)+\left(
\begin{array}{c}
C_0e_+ \\
iC_0e_-z +C_1\\
\end{array}\right)\right], \quad z\in {\Bbb C}\setminus[-1,1].
\label{7.16}
\eeq
where $C_0=-P_*/(2\pi a)$, and $\BGY(z)$ is simplified to the form
\beq
\BGY(z)=-\fr{1}{2\pi(\nu+1)e_0}\int_{-1}^1\fr{\BT(\tau)g(\tau)d\tau}{\GL^+(\tau)(\tau-z)}.
\label{7.17}
\eeq
It is seen that the  first component of the vector $\BGF(z)$,
the function $\GF_1(z)$, has a nonintegrable singularity of order $-5/4$  at the point $z=1$ unless 
\beq
C_1=-\GY_2(1)+\fr{ie_-P_*}{2\pi a}.
\label{7.18}
\eeq
This condition is necessary and sufficient for the function $\GF_1(t)$ being integrable in the vicinity of the point $z=1$.
Notice that the condition (\ref{7.18}) coincides with  (\ref{4.26}) when $k=1$. Now taking the limit $k\to 1$
in the representation formula (\ref{5.18}) for the solution of the integral equation (\ref{3.20})
we obtain
$$
\psi(at)=\fr{i(\nu-1)g(t)}{\nu_1}-\fr{2i}{(1+\nu)\sqrt{\nu_0}}(1-t)^{-3/4-i\Gg}(1+t)^{-1/4+i\Gg}
$$
\beq
\times\left[\GY_1(t)+\fr{\GY_2(t)-\GY_2(1)}{\sqrt{1-t^2}}-
\fr{P_*}{2\pi a}\left(e_+-ie_-\sqrt{\fr{1-t}{1+t}}\right)\right], \quad -1<t<1.
\label{7.19}
\eeq
Utilizing formula (\ref{7.17}) we find the principal values of the integrals $\Psi_1(t)$ and
$(1-t^2)^{-1/2}[\Psi_2(t)-\Psi_2(1)]$. They read
$$
\Psi_1(t)=\fr{1}{2\pi (\nu+1)e_0^2}\int_{-1}^1\fr{\Go_1(\tau)g(\tau)d\tau}{\tau-t}, 
$$
\beq
\fr{\Psi_2(t)-\Psi_2(1)}{\sqrt{1-t^2}}=\fr{1}{2\pi (\nu+1)e_0^2}\sqrt{\fr{1-t}{1+t}}
\int_{-1}^1\fr{\Go_2(\tau)g(\tau)d\tau}{\tau-t},  \quad -1<t<1,
\label{7.20}
\eeq
where
\beq
\Go_1(t)=(1-t)^{3/4+i\Gg}(1+t)^{1/4-i\Gg},\quad 
\Go_2(t)=(1-t)^{1/4+i\Gg}(1+t)^{3/4-i\Gg}.
\label{7.21}
\eeq
We next substitute formulas (\ref{7.20}) into (\ref{7.19}). The result,
$$
\psi(at)=\fr{i(\nu-1)g(t)}{\nu_1}-
\fr{1}{\pi\Go_1(t)}\left[\fr{1}{\nu_1}\int_{-1}^1\fr{\Go_1(\tau)g(\tau)d\tau}{\tau-t}-\fr{iP_*e_+}{a(1+\nu)\sqrt{\nu_0}}\right]
$$
\beq
- \fr{1}{\pi\Go_2(t)}\left[\fr{1}{\nu_1}\int_{-1}^1\fr{\Go_2(\tau)g(\tau)d\tau}{\tau-t}-\fr{P_*e_-}{a(1+\nu)\sqrt{\nu_0}}\right],
\label{7.22}
\eeq
is the exact solution to the integral equation (\ref{7.8}), the limit of the solution (\ref{5.18}) as $k\to 1$.

\subsection{The limiting case $k=1$ when $w^\circ=0$, $\Gs^\circ=0$, and $h=\const$} 

Set now $w^\circ(x)=0$, $h'(x)=0$, $-b<x<b$, and  $\Gs^\circ(x)=0$, $-a<x<a$. This case was treated in  {\bf(\ref{mus})}  by the method of complex potentials.
We aim to simplify the representation formula for the solution (\ref{7.22}) derived for the limiting case $k=1$. 
We have $P_*=P$, and formula (\ref{7.22}) reads
$$
\psi(at)=-\fr{P}{2\pi a}\left[
\fr{i(1-\nu)^2}{\nu_1\sqrt{1-t^2}}-\fr{2}{(1+\nu)\sqrt{\nu_0}}
\left(\fr{ie_+}{\Go_1(\tau)}+\fr{e_-}{\Go_2(\tau)}\right)
\right.
$$
\beq
\left.
+\fr{1-\nu}{\pi\nu_1}
\left(\fr{1}{\Go_1(t)}
\int_{-1}^1
\fr{\Go_1(\tau)d\tau}{\sqrt{1-\tau^2}(\tau-t)}+\fr{1}{\Go_2(t)}
\int_{-1}^1
\fr{\Go_2(\tau)d\tau}{\sqrt{1-\tau^2}(\tau-t)}\right)
\right].
\label{7.24}
\eeq
The two integrals in (\ref{7.24}) can be evaluated by
means of the relation
\beq
\int_{-1}^1\fr{(1-\tau)^\Ga(1+\tau)^{-\Ga}d\tau}{\tau-t}
=\pi\cot\pi\Ga\left(\fr{1-t}{1+t}
\right)^\Ga-\fr{\pi}{\sin\pi\Ga}, \quad -1<t<1,
\label{7.25}
\eeq
a particular case of the more general formula  {\bf(\ref{tri})}
$$
\fr{1}{\pi}\int_{-1}^1 (1-\tau)^\Ga(1+\tau)^\Gb P_n^{(\Ga,\Gb)}(\tau)\fr{d\tau}{\tau-t}
=\cot\pi\Ga (1-t)^\Ga(1+t)^\Gb P_n^{(\Ga,\Gb)}(t)
$$
$$
-\fr{2^{\Ga+\Gb}\GG(\Ga)\GG(n+\Gb+1)}{\pi\GG(n+\Ga+\Gb+1)}
F\left(n+1,-n-\Ga-\Gb; 1-\Ga; \fr{1-t}{2}\right),
$$
\beq
\Ga>-1, \quad \Gb>-1; \quad \Ga\ne 0,1,2,\ldots; \quad -1<t<1,
\label{7.25'}
\eeq
obtained by employing the integral representation of the Jacobi functions of the second kind
in terms of the Jacobi polynomials $ P_n^{(\Ga,\Gb)}(t)$. Here, $F$ is the hypergeometric function.
Alternatively, formula (\ref{7.25'}) can be derived by converting the left hand-side into an integral
of the Mellin convolution type and applying the theory of residues  {\bf(\ref{pop})}.

On plugging the expressions obtained into formula (\ref{7.24}) and using the identities
$$
\fr{e_-}{e_+}+\fr{e_+}{e_-}=1-\nu,
$$
\beq
\fr{\nu-1}{\nu_1 e_\pm}-\fr{2i e_\pm}{(1+\nu)\sqrt{\nu_0}}
=\pm\fr{1}{(1+\nu)e_0^2 e_\pm} 
\label{7.26}
\eeq
 we deduce
\beq
\psi(at)=-\fr{P}{2\pi a(1+\nu)\sqrt{\nu_0}}\left[\fr{1}{e_-} (1-t)^{-1/4-i\Gg}(1+t)^{-3/4+i\Gg}
-\fr{i}{e_+} (1-t)^{-3/4-i\Gg}(1+t)^{-1/4+i\Gg}\right],
\label{7.27}
\eeq
where the parameters $e_\pm$ and $\Gg$ are given by
(\ref{7.14}) and (\ref{4.13}), respectively.
It is directly verified that this function satisfies the equilibrium condition (\ref{3.22}) ($k=1$).
To show that the function (\ref{7.27}) is a solution of equation (\ref{7.8}), in  addition to the relation (\ref{7.25}), we employ the formula 
\beq
\int_{-1}^1\fr{(1-\tau)^\Ga(1+\tau)^{-\Ga-1}d\tau}{\tau-t}
=\pi\cot\pi\Ga(1-t)^\Ga(1+t)^{-\Ga-1}, \quad -1<t<1,
\label{7.28}
\eeq
that can be deduced from (\ref{7.25'}). Then, upon substitution, one can check that the function  (\ref{7.27})
solves the integral equation (\ref{7.8}) with the function $g(t)$ given by (\ref{6.11}).

\setcounter{equation}{0}
  
\section{Conclusions}

In this work we have analyzed two model contact problems on  a rigid inclusion debonded from an
elastic medium. Model 1 concerns an inclusion completely debonded from the matrix, and the crack
formed in the upper side of the inclusion 
penetrates into the medium. In Model 2, the crack length $2a$ is less than the inclusion length $2b$.
Each model is governed by a singular integral equation with the same kernel but a different
right-hand side.
We have developed a method that ultimately leads to a closed-form solution of the integral equation.
The main feature of the method is the solution of the associated order-2 vector Riemann-Hilbert problem
with the Chebotarev-Khrapkov matrix coefficient in a finite segment, not in a closed contour
(an infinite line) as the classical Khrapkov scheme {\bf(\ref{khr})} requires. We have examined the behavior 
of the solution at the crack and inclusion tips, determined the SIFs and the contact stresses
and reported sample numerical results for them. To verify the numerical results, we obtained 
the SIFs for the limiting case $k=b/a\to 0$ independently. It turns out that  both factors, $K_I$ and $K_{II}$,
in the case $k=0^+$ are very close to the values computed from the general formulas for $k$ small ($k=0.005$).
   
By passing to the limit $k\to 1$ in the solution found we have managed to derive a closed-form solution 
for the particular case $b=a$. That solution coincides with the one known in the literature {\bf(\ref{mus})}.
It turns out that the traction components have  square root 
singularities at the crack and inclusion tips in the cases $a<b$ and $a>b$, while 
they have a stronger singularity of order  $-3/4$ in the limiting case $a=b$. Also, when $a\to b^\pm$,
the absolute values of the  Mode-I and II SIFs grow to infinity. This fact allows us to conclude that 
the case $a=b$ is unstable: when the  crack formed on the upper surface of the inclusion ($a<b$)
starts growing and its tips approach the inclusion ends ($a\to b$ and $a<b$),
the SIFs grow unboundly, and the crack tends to speed up to pass the inclusion tips and  penetrate into the matrix.
After that moment the SIFs decrease, and eventually the crack  stops
at some distance from the inclusion tips.

\vspace{.2in}

\newpage

{\centerline{\Large\bf  References}}

\vspace{.1in}

\begin{enumerate}

\item\label{she} D. I. Schermann,  Probl\`eme mixte de la th\'eorie du potentiel et de la th\'eorie de l'\'elasticit\'e 
pour un plan ayant un nombre fini de coupures rectilignes, {\it C. R. (Doklady) Acad. Sci. URSS (N.S.)} 
{\bf 27} (1940) 329-333. 

\item\label{mus} N. I. Muskhelishvili,  {\it Some basic problems of the mathematical theory of elasticity}
(P. Noordhoff, Ltd., Groningen 1963).

\item\label{cher} G. P. Cherepanov, 
Solution of a linear boundary value problem of Riemann for two functions and its application to certain mixed problems in the plane theory of elasticity, {\it J. Appl. Math. Mech.} {\bf 26} (1962) 1369-1377.

\item\label{pop} G. Ya. Popov, {\it Concentration of elastic stresses near stamps, cuts, thin inclusions and  reinforcements} (Nauka, Moscow, 1982).

\item\label{zve} E. I. Zverovich,  A mixed problem of elasticity theory for the plane with cuts that lie on the real axis. 
In {\it Proc.  Symp. on Continuum Mechanics and Related Problems of Analysis (Tbilisi, 1971)} {\bf 1}, 
103-114 (Mecniereba, Tbilisi 1973).

\item\label{sil} V. V. Silvestrov, The method of Riemann surfaces in the problem of interface cracks and inclusions under concentrated forces, {\it Russian Math. (Iz. VUZ)}  {\bf 48} (2005) no. 7, 75-88.

\item\label{ant} Y.A. Antipov, A delaminated inclusion in the case of adhesion and slippage, {\it J. Appl. Math. Mech.} {\bf 60} (1996) 665-675.

\item\label{bar} D. I. Bardzokas and S. M. Mkhitaryan, On interaction of various types of stress concentrators. In {\it Proc. Int. Conf. on Problems in Continuum Mechanics}, 6 p. (Erevan, 2007). 

\item\label{sel} A. P. S. Selvadurai and B. M. Singh,
On the expansion of a penny-shaped crack by a rigid circular disc
inclusion, {\it Int. J.  Fracture} {\bf 25} (1984) 69-77.

\item\label{hwu}  C. Hwu, Y. K. Liang and  W. J. Yen,
 Interactions between inclusions and various types of cracks,
  {\it Int. J.  Fracture} {\bf 73} (1995) 301-323.

\item\label{che} G. N. Chebotarev, On closed-form solution of a Riemann boundary value problem for $n$ pairs
of functions, {\it Uchen. Zap. Kazan. Univ.} {\bf 116} (1956) 31-58.

\item\label{khr} A. A. Khrapkov, Certain cases of the elastic equilibrium of an infinite wedge with a non- symmetric notch at the vertex, subjected to concentrated forces, {\it J. Appl. Math. Mech.} {\bf 35} (1971) 625-637.

\item\label{gak} F. D. Gakhov, {\it Boundary Value Problems} (Pergamon Press, Oxford, 1966).

\item\label{cher2} G. P. Cherepanov, {\it Mechanics of brittle fracture}  (McGraw-Hill, New York, 1979).

\item\label{tri} F. G. Tricomi, On the finite Hilbert transformation, {\it Quart. J. Math.} {\bf 2} (1951), 199-211.

\end{enumerate}

\vspace{.2in}

\end{document}